\documentclass[12pt,leqno,twoside]{article}
\usepackage{amssymb,amsmath}
\usepackage{comment}
\usepackage{color}
\usepackage[dvipdfmx]{graphicx}
\usepackage[mathscr]{eucal}
\usepackage{fancyhdr}
\newtheorem{df}{Definition}[section]
\newtheorem{thm}[df]{Theorem} 
\newtheorem{cor}[df]{Corollary} 
\newtheorem{prop}[df]{Proposition} 
\newtheorem{lem}[df]{Lemma} 
\newtheorem{rem}[df]{Remark} 
\makeatletter
 
\@addtoreset{equation}{section} 

\let\@fnsymbol\@arabic

\makeatother
\newcommand{\qed}{$\hfill\square$}
\newcommand{\dis}{\displaystyle}
\renewcommand{\phi}{\varphi}

\newcommand{\D}{\Delta}
\newcommand{\G}{\Gamma}
\newcommand{\bN}{{\mathbb{N}}}
\newcommand{\bR}{{\mathbb{R}}}

\newcommand{\bme}{\mbox{\boldmath$e$}}

\newcommand{\bmx}{\mbox{\boldmath$x$}}
\newcommand{\bmy}{\mbox{\boldmath$y$}}
\newcommand{\bmz}{\mbox{\boldmath$z$}}
\newcommand{\bmxi}{\mbox{\boldmath$\xi$}}
\newcommand{\bmzeta}{\mbox{\boldmath$\zeta$}}
\newcommand{\bmeta}{\mbox{\boldmath$\eta$}}

\newcommand{\bmtau}{\mbox{\boldmath$\tau$}}

\newcommand{\bmn}{\mbox{\boldmath$n$}}

\newcommand{\bmzero}{\mbox{\boldmath$0$}}
\newcommand{\sbme}{\mbox{{\scriptsize\boldmath$e$}}}
\newcommand{\sbmx}{\mbox{{\scriptsize\boldmath$x$}}}
\newcommand{\sbmy}{\mbox{{\scriptsize\boldmath$y$}}}
\newcommand{\sbmz}{\mbox{{\scriptsize\boldmath$z$}}}
\newcommand{\sbmzeta}{\mbox{{\scriptsize\boldmath$\zeta$}}}
\newcommand{\sbmeta}{\mbox{{\scriptsize\boldmath$\eta$}}}
\newcommand{\sbmxi}{\mbox{{\scriptsize\boldmath$\xi$}}}

\newcommand{\sbmzero}{\mbox{{\scriptsize\boldmath$0$}}}
\newcommand{\tbmx}{\mbox{{\tiny\boldmath$x$}}}
\newcommand{\sW}{{\mathscr{W}}}

\newcommand{\al}{\alpha}
\newcommand{\be}{\beta}
\newcommand{\om}{\Omega}

\newcommand{\p}{\partial}

\newcommand{\bC}{\mathbb{C}}

\newcommand{\bZ}{\mathbb{Z}}

\newcommand{\ms}{\hspace*{-1.5pt}}
\renewcommand{\,}{\hspace*{1pt}}

\allowdisplaybreaks[1]

\begin{document}

\title{A thresholding algorithm to Willmore-type flows \\[0.1cm] 
via fourth order linear parabolic equation}


\author{Katsuyuki Ishii, 
Yoshihito Kohsaka, 
Nobuhito Miyake \\[0.05cm]
and Koya Sakakibara}

\date{}

\maketitle

\pagestyle{fancy}
\fancyhead{}
\fancyhead[LE]{\thepage}
\fancyhead[CE]{{\sc K. Ishii, Y. Kohsaka, N. Miyake and K. Sakakibara}}
\fancyhead[RO]{\thepage}
\fancyhead[CO]{{\sc A thresholding algorithm to Willmore-type flows}}
\renewcommand{\headrulewidth}{0.0pt}
\fancyfoot{}

\vspace*{-0.325cm}

\begin{abstract}
We propose a thresholding algorithm to Willmore-type flows in $\mathbb{R}^N$. 
This algorithm is constructed based on the asymptotic expansion of the solution to the initial 
value problem for a fourth order linear parabolic partial differential equation whose initial data is the indicator 
function on the compact set $\Omega_0$. The main results of this paper demonstrate 
that the boundary $\partial\Omega(t)$ of the new set $\Omega(t)$, generated by our algorithm, is included in 
$O(t)$-neighborhood of $\partial\Omega_0$ for small $t>0$ and that the normal velocity from $ \partial\Omega_0 $ 
to $ \partial\Omega(t) $ is nearly equal to the $L^2$-gradient of Willmore-type energy 
for small $ t>0 $.  
Finally, numerical examples of planar curves governed by the Willmore flow are provided by using 
our thresholding algorithm. 
\end{abstract}

\vspace*{-0.75cm}

\footnote[0]{This project is supported by JSPS KAKENHI Grant Number JP23H00085, JP23K03215, JP19K03562 (for K.\,I. and Y.\,K.), 
JP22KJ0719 (for N.\,M.) and JP22K03425 (for K.\,S.).}


\section{Introduction}

In this paper we propose a thresholding algorithm to the $L^2$-gradient flow of Willmore-type energy 
in $\mathbb{R}^N$. 
Let $\Gamma$ be a hypersurface in $\mathbb{R}^N$. The Willmore-type energy $\mathcal{E}_{\lambda}^N(\Gamma)$ 
for $\Gamma$ is defined as 
\[
\mathcal{E}_{\lambda}^N(\Gamma)
=\left\{\begin{array}{ll}
\displaystyle\frac12\int_\Gamma \kappa^2\,ds+\lambda\int_\Gamma\,ds&\text{if}\,\ N=2, \\[0.5cm]
\displaystyle\frac12\int_\Gamma H^2\,dS_\Gamma-\frac23\int_\Gamma\sum_{\substack{i,j\in\Lambda \\ i<j}}
\kappa_i\kappa_j\,dS_\Gamma
+\lambda\int_\Gamma\,dS_\Gamma&\text{if}\,\ N\ge3,
\end{array}\right.
\]
where $\lambda\in\mathbb{R}$. 
For $N=2$, $\kappa$ is the curvature of a planar curve $\Gamma$ and $s$ is the arc-length parameter. 
For $N\ge3$, $ \kappa_i $\,($ i\in\Lambda:=\{1,2,\cdots,N-1\} $) are the principal curvatures of $ \G $ and 
$ H $ is the ($ (N-1) $-times) mean curvature of $ \G $. 
This can be regarded as a generalization of the Willmore energy. 
Note that the energy $ \mathcal{E}_0^N(\G) $ with $ N\geq 3 $ appears in the asymptotic expansion of the heat content 
(see Angiuli--Massari--Miranda\,\cite{An;Ma;Mi;13}).  

Let $ \{\G(t)\}_{t\geq 0} $ be a family of hypersurfaces in $ \bR^N $ and assume that the motion of $ \G(t) $ 
is governed by the $L^2$-gradient flow
\begin{equation}\label{geometric-flow}
V=-\nabla_{L^2}\mathcal{E}_{\lambda}^N(\Gamma(t)),
\end{equation}
where $ V $ is the normal velocity of $ \G(t) $ and $\nabla_{L^2}\mathcal{E}_{\lambda}^N(\Gamma)$ is the $L^2$-gradient 
of $\mathcal{E}_{\lambda}^N(\Gamma)$ given by
\begin{equation}\label{L2-gradient}
\nabla_{L^2}\mathcal{E}_{\lambda}^N(\Gamma)
=\left\{\begin{array}{ll}
\p^2_s\kappa+\dfrac{1}{2}\kappa^3-\lambda\kappa&\mbox{if}\ N=2, \\[0.325cm]
\D_gH+H|A|^2-\dfrac{1}{2}H^3-\lambda H&\mbox{if}\ N=3, \\[0.325cm]
\dis\D_gH+H|A|^2-\dfrac{1}{2}H^3+2\sum_{\substack{i,j,k\in\Lambda \\ i<j<k}}\kappa_i\kappa_j\kappa_k-\lambda H
&\mbox{if}\ N\ge4.
\end{array}\right.  
\end{equation}
Here $ \p^2_s $ is the second order differential operator with respect to $s$, 
$ \D_g $ is the Laplace-Beltrami operator on $ \G $ by the induced metric $g=(g_{ij})$ and $ |A|^2 $ denotes the norm 
of the second fundamental form $ A=(h_{ij}) $, 
which is defined as $ |A|^2:=\sum\limits_{i,j,k,\ell\in\Lambda} g^{ij}g^{k\ell}h_{ik}h_{j\ell}
\left(=\sum\limits_{i\in\Lambda}\kappa_i^2\right) $ where $g^{-1}=(g^{ij})$ is the inverse matrix of $g=(g_{ij})$. 
If $ \G(t) $ is embedded and encloses a domain $ D(t) $, we choose the orientation induced by the outer unit normal, 
so that $ V $ is positive if $ D(t) $ grows and $ H $ is negative if $\Gamma(t)$ is a spherical surface.  
Note that the term $\sum\limits_{i,j,k\in\Lambda}\kappa_i\kappa_j\kappa_k $ for $ N\ge4 $ in \eqref{L2-gradient} is 
derived as the first variation of the integral of $ (-1/3)\sum\limits_{i,j\in\Lambda}\kappa_i\kappa_j $ on $\Gamma$ 
(see e.g. Reilly\,\cite{rei;73}). Also note that if $N=3$ and the topology of $\Gamma$ is fixed, 
the integral of $\kappa_1\kappa_2$ on $\Gamma$ is constant by virtue of the Gauss--Bonnet theorem, so that 
its first variation is zero. 


The equation \eqref{geometric-flow} with $ N=2,3 $ and $\lambda=0$ is the Willmore flow (WF for short). 
For the results of the existence and the asymptotic behavior of the WF 
and related flows, see e.g. Simonett\,\cite{sim;01}, Kuwert--Sch\"{a}tzle\,\cite{kuw;sch;01,kuw;sch;02}, 
Dziuk--Kuwert--Sch\"{a}tzle\,\cite{dzi;kuw;sch;02}, Okabe--Wheeler \cite{oka;whe;21}, Rupp\,\cite{rup;21} 
and references therein.  As for the approximation schemes and the methods of numerical computations of the flow 
by \eqref{geometric-flow}, there are many results taking account of various applications.  
Mayer--Simonett\,\cite{may;sim;02} is one of the first numerical approaches for the WF in $ \bR^3 $.  They used 
a finite difference scheme to the WF and numerically observed that the WF can develop singularities 
in finite time. Rusu\,\cite{rus;05} presented an algorithm to the WF in $ \bR^3 $ based on the variational method and  
studied a semi- and a fully discrete scheme in space and a semi-implicit method in time. 
Dziuk\,\cite{dzi;08} introduced a parametric finite element method to the WF in general space dimensions. 
In \cite{bar;gar;nur;07,bar;gar;nur;08} etc., Barrett, Garcke and N\"{u}rnberg studied parametric 
finite element methods for fourth order geometric evolution problems such as surface diffusion flow and the WF.  
Furthermore, it is well known that the WF can be approximated by the fourth order phase-field equations 
or equivalent systems of PDE's, which are derived from approximations of the Willmore functional. 
Loreti--March\,\cite{lor;mar;00} obtained 
the formal asymptotic expansions of solutions of the fourth order phase-field equations (or equivalent 
systems) and derived \eqref{geometric-flow} for $N=3$.  Bretin--Masnou--Oudet\,\cite{bre;mas;oud;15} gave  
similar results to some modified versions of the fourth order phase-field equations and the related energies.   
They also presented in \cite{bre;mas;oud;15} some numerical simulations of the flows for $ N=2,3 $, based on 
their formal asymptotic expansions.  Colli--Lauren\c{c}ot\,\cite{col;lau;11,col;lau;12} studied the well-posedness 
of a phase-field approximation to the WF with volume and area constraints in $\mathbb{R}^N\,(1\le N\le3)$. 
R\"{a}tz--R\"{o}ger\,\cite{rat;rog;21} introduced a new diffuse-interface approximation of the WF avoiding intersections 
of phase boundaries that do not correspond to the intended sharp interface evolution. They also justified 
the approximation property by a Gamma convergence for the energies and a matched asymptotic expansion 
for the flow. Fei--Liu\,\cite{fei;liu;21} rigorously proved the convergence 
of the zero level set of the solutions of the phase-field system to the WF for $ N=2,3 $ if the smooth WF exists.   

The purpose of this paper is to introduce a thresholding algorithm by using the following Cauchy problem for 
the fourth order linear parabolic equation:
\begin{equation}\label{4th_diffusion}
\left\{\begin{array}{l}
u_t=-\Delta^2u+\lambda\Delta u\,\,\ \text{in}\ \mathbb{R}^N\times(0,\infty), \\
u(\bmx,0)=\chi_{\Omega_0}(\bmx):=\left\{\begin{array}{l} 1\,\,\ 
\text{in}\ \Omega_0, \\ 0\,\,\ \text{in}\ \mathbb{R}^N\setminus \Omega_0. \end{array}\right.
\end{array}\right.
\end{equation}
Here $N\ge2$, $\lambda\in\bR$, and $\om_0\subset\bR^N$ is compact set with smooth boundary. 
By the derivation of a threshold function from the solution to the above problem, we obtain a thresholding algorithm 
to the motion of $\G(t)$ by \eqref{geometric-flow} at least formally.  
The outline of our algorithm is as follows: 
Set $ h>0 $ as a time step. For a given compact set $\Omega_0$ in $\mathbb{R}^N$ with a smooth boundary 
$\partial \Omega_0$, we solve the initial value problem \eqref{4th_diffusion}.
Next, let $ u^0 $ be the solution to \eqref{4th_diffusion} and set $ u^0_a(\bmx,t):=u^0(\bmx,a^4t) $ for $ a>0 $. 
Then, define a threshold function $U^0$ as
\[
     U^0(\bmx,t):=u^0_{3a}(\bmx,t)-3u^0_{2a}(\bmx,t)+3u^0_a(\bmx,t)
\] 
and give a new set $ \om_1 $ by 
\[
     \om_1:=\left\{\bmx\in\bR^N\;|\;U^0(\bmx,h)\geq\frac{1}{2}\right\}.  
\]
Repeating this procedure inductively, we obtain a sequence $ \{\om_k\}_{k=0}^{\infty} $ 
of compact subsets of $ \bR^N $.   Set 
\[
    \om^h(t):=\om_k \quad\mbox{for}\,\ kh\leq t<(k+1)h,\,\ k=0,1,2,\cdots.
\]
Then, letting $ h\to 0 $, we observe at least formally that $ \om^h(t) $ converges to a compact set $ \om(t)(\subset\bR^N) $ 
and that $ \p\om(t) $ moves by \eqref{geometric-flow} if we choose a suitable constant $ a $. 
In order to justify the thresholding algorithm explained above, we derive the asymptotic expansion of the solution to 
\eqref{4th_diffusion} near $ \p\om_0 $.
For the details of the justification, see the argument in subsection \ref{subsec:proposed_algorithm}.

Thresholding algorithms to the geometric evolution equations was firstly introduced by 
Bence--Merriman--Osher\,\cite{ben;mer;osh;92} to numerically compute the mean curvature flows.  
Based on the level set approach for geometric evolution equations, the convergence and generalizations of their algorithm 
were studied by Mascarenhas\,\cite{mas;92}, Evans\,\cite{eva;93}, Barles--Georgelin\,\cite{ba;ge;95}, 
Ishii\,\cite{is;95}, Ishii--Pires--Souganidis\,\cite{is;pi;sou;99}, Vivier\,\cite{vi;00}, Leoni\,\cite{le;01} and so on. 
Recently, another approach was suggested by Esedo\=glu--Otto\,\cite{ese;ott;15}, which was considered 
the thresholding algorithm for the multi-phase mean curvature flow. 
They gave the interpretation such that the thresholding algorithm can be regarded as a minimizing movement scheme. 
For the result of the convergence and the further development on this approach, 
see Laux--Otto\,\cite{lau;ott;16,lau;ott;20-1,lau;ott;20-2}, Laux--Lelmi\,\cite{lau;lel;22} and 
Fuchs--Laux\,\cite{fuc;lau;24}.

On the thresholding algorithm to the WF, there are results by Grzhibovskis--Heintz\,\cite{grz;hei;08} in $\mathbb{R}^3$ 
and Esedo\={g}lu--Ruuth--Tsai\,\cite{ese;ruu;tsa;08} in $\mathbb{R}^2$. 
In \cite{grz;hei;08} and \cite{ese;ruu;tsa;08} the asymptotic expansion of the convolution 
$ (t^{-N/4}\rho(|\cdot|/t^{1/4})\ast\chi_{\om_0})(\bmx) $ is used to define a threshold function. 
Here $t^{-N/4}\rho(|\cdot|/t^{1/4})$ is a modified Gauss kernel or some similar ones. Note that in \cite{ese;ruu;tsa;08} 
the $L^2$-gradient flow of the Helfrich functional in $\mathbb{R}^2$ was also considered. 
On the details of the difference between their thresholding algorithm and ours, 
see Remark \ref{rem:difference}.
Metivet--Sengers--Isma\"{i}l--Maitre\,\cite{met;sen;ism;mai;21} treated the diffusion-redistance scheme in 
$\mathbb{R}^2$ or $\mathbb{R}^3$, which is a variant of the algorithm by \cite{ese;ruu;tsa;08}. 

Referring to \cite{grz;hei;08} and \cite{ese;ruu;tsa;08}, the space-time scale $|\bmx|/t^{1/4}$ 
plays a key role to obtain the WF from the formal asymptotic expansions of their convolutions.  
Indeed, in \cite{grz;hei;08} and \cite{ese;ruu;tsa;08}, they used a modified Gauss kernel whose 
space-time scale is $|\bmx|/t^{1/4}$ instead of a natural space-time scale $|\bmx|/t^{1/2}$ 
of the usual Gauss kernel (see Remark \ref{rem:difference} on the details of the calculation). 
Based on this fact, we arrive at the idea of using the fundamental solution to the fourth order linear parabolic equation 
in \eqref{4th_diffusion} to construct a thresholding algorithm to the flow by \eqref{geometric-flow} 
since a natural space-time scale of its solution is $|\bmx|/t^{1/4}$. 
Another reason to consider the equation in \eqref{4th_diffusion} is that \eqref{geometric-flow} 
can be written by the signed distance function $d$ as follows 
(cf. \cite[Section 3.3]{bre;mas;oud;15} and \cite[Lemma A.2]{fei;liu;21}):  
\[
d_t=-\D^2d+\langle D(\langle D(\D d),Dd\rangle), Dd\rangle+2\langle D(\D d),Dd\rangle\D d
+\frac{1}{2}(\D d)^3+\lambda \D d 
\,\,\ \mbox{on}\,\ \G(t).  
\]
Dropping all of the nonlinear terms, we have the fourth order linear parabolic equation 
in \eqref{4th_diffusion} and hence it is regarded as the "rough" approximation of this equation. 
Furthermore, as the benefit of using the fundamental solution to \eqref{4th_diffusion}, 
the part related to the area constraint is naturally derived from the term of Laplacian of the equation 
in \eqref{4th_diffusion}. Such a derivation considering the structure is difficult if we use a modified 
Gauss kernel. We remark that for a thresholding algorithm to the WF, at present, there are no results 
on the convergence to some suitable solution to the WF. Since it seems that our approach is 
more natural compared with that in \cite{grz;hei;08} and \cite{ese;ruu;tsa;08}, it is expected that 
the construction of some suitable solution and the convergence to it are shown based on our results. 
Indeed, in order to prove the convergence, the ideas based on the gradient flow as in 
\cite{fuc;lau;24,lau;lel;22,lau;ott;16,lau;ott;20-1,lau;ott;20-2} may be useful for the WF. 
We think that our algorithm has a strong possibility for connecting to their approach.

This paper is organized in the following way.  In section 2 we derive some formulae and pointwise 
estimates of the fundamental solution $ G_{N,\lambda} $ to the operator $ \p_t+(-\D)^2+\lambda(-\D)$ 
and of its derivatives.  
In section 3 we discuss the formal asymptotic expansion of the solution to \eqref{4th_diffusion}, 
which is stated in Theorem \ref{thm:asympt_expansion1}.
Section 4 is devoted to the justification of our algorithm. 
In subsection \ref{subsec:proposed_algorithm} we recall the algorithm and the fact to be justified 
and prove this fact in subsection \ref{sec:evolving-sets}. 
Section \ref{sec:numerical_experiments} presents the results of numerical experiments. 
We consider the equation in \eqref{4th_diffusion} with $\lambda=0$ and $N=2$, i.e. the case corresponding to 
the WF in $\mathbb{R}^2$, in a periodic square domain. 
The final section is Appendix.  


\section{Preliminaries}

In this section we derive some properties and estimates of the fundamental solution $G_{N,\lambda}$ to 
the operator $\partial_t+(-\Delta)^2+\lambda(-\Delta)$ on $\mathbb{R}^N\times(0,\infty)$ where 
$\Delta$ is the Laplace operator on $\mathbb{R}^N$ and $\lambda\in\mathbb{R}$. 

Define the Fourier transform as 
\[
\mathcal{F}_N[\psi](\bmxi):=\int_{\mathbb{R}^N}\psi(\bmx)e^{-\textbf{i}\langle\sbmxi,\sbmx\rangle_N}\,d\bmx, \quad
\mathcal{F}_N^{-1}[\psi](\bmx):=c_N\int_{\mathbb{R}^N}\psi(\bmxi)e^{\textbf{i}\langle\sbmx,\sbmxi\rangle_N}\,d\bmxi
\]
for $ \psi\in L^1(\bR^N) $ where $\textbf{i}:=\sqrt{-1}$, $c_N:=(2\pi)^{-N}$ and 
$\langle\,\cdot\,,\,\cdot\,\rangle_N$ is the inner product on $\mathbb{R}^N$. Then $G_{N,\lambda}$ is given by 
\[
G_{N,\lambda}(\bmx,t)=\mathcal{F}_N^{-1}[e^{-(|\,\cdot\,|^4+\lambda|\,\cdot\,|^2)t}](\bmx)
=c_N\int_{\mathbb{R}^N}e^{-(|\sbmxi|^4+\lambda|\sbmxi|^2)t+\textbf{i}\langle\sbmx,\sbmxi\rangle_N}\,d\bmxi, 
\]
and we readily see that
\begin{equation}\label{integralG=1}
\int_{\mathbb{R}^N}G_{N,\lambda}(\bmx,t)\,d\bmx
=\mathcal{F}_N[\mathcal{F}_N^{-1}[e^{-(|\,\cdot\,|^4+\lambda|\,\cdot\,|^2)t}]](\bmzero)
=e^{-(|\sbmxi|^4+\lambda|\sbmxi|^2)t}\big|_{\sbmxi=\sbmzero}=1.
\end{equation}
%
%
Set $g_N(\bmx):=G_{N,0}(\bmx,1)$. That is,  
\begin{equation}\label{g_N}
g_N(\bmx)=\mathcal{F}_N^{-1}[e^{-|\,\cdot\,|^4}](\bmx)
=c_N\int_{\mathbb{R}^N}e^{-|\sbmxi|^4+\textbf{i}\langle\sbmx,\sbmxi\rangle_N}\,d\bmxi. 
\end{equation}
We derive the expansion of $ G_{N,\lambda} $ by use of $ g_N $ and its derivatives. 

\begin{prop}\label{prop:rep_4th-laplace}
$G_{N,\lambda}$ is represented as
\[
G_{N,\lambda}(\bmx,t)=
\frac1{t^{N/4}}\sum_{m=0}^\infty\frac{(-\lambda)^mt^{m/2}}{m!}(-\Delta_{\sbmz})^mg_N\Bigl(\frac{\bmx}{t^{1/4}}\Bigr),
\]
where $\Delta_{\sbmz}\,g_N$ is the Laplacian of $g_N(\bmz)$ with respect to $ \bmz\in\bR^N $. 
\end{prop}
\noindent
{\bf Proof.} 
Applying the change of variable $\bmxi=\bmeta/t^{1/4}$, we have
\begin{align*}
G_{N,\lambda}(\bmx,t)
=&\,\frac{c_N}{t^{N/4}}\int_{\mathbb{R}^N}e^{-|\sbmeta|^4+\textbf{i}\langle\frac{\sbmx}{t^{1/4}},\sbmeta\rangle_N}
e^{-\lambda|\sbmeta|^2t^{1/2}}\,d\bmeta \\
=&\,\frac{c_N}{t^{N/4}}\int_{\mathbb{R}^N}\sum_{m=0}^\infty\frac{(-\lambda)^mt^{m/2}}{m!}|\bmeta|^{2m}
e^{-|\sbmeta|^4+\textbf{i}\langle\frac{\sbmx}{t^{1/4}},\sbmeta\rangle_N}\,d\bmeta. 
\end{align*}
Since
\begin{align*}
\int_{\mathbb{R}^N}\sum_{m=0}^\infty\left|\frac{(-\lambda)^mt^{m/2}}{m!}|\bmeta|^{2m}e^{-|\sbmeta|^4+\textbf{i}\langle
\frac{\sbmx}{t^{1/4}},\sbmeta\rangle_N}\right|d\bmeta
&=\int_{\mathbb{R}^N}\sum_{m=0}^\infty\frac{|\lambda|^mt^{m/2}}{m!}|\bmeta|^{2m}e^{-|\sbmeta|^4}d\bmeta \\
&=\int_{\mathbb{R}^N}e^{-|\sbmeta|^4+|\lambda||\sbmeta|^2t^{1/2}}d\bmeta<+\infty\nonumber,
\end{align*}
the Lebesgue convergence theorem implies that the termwise integration is possible and hence
\begin{equation}\label{rep_G_Nlambda}
G_{N,\lambda}(\bmx,t)=\frac{c_N}{t^{N/4}}\sum_{m=0}^\infty\frac{(-\lambda)^mt^{m/2}}{m!}
\int_{\mathbb{R}^N}|\bmeta|^{2m}e^{-|\sbmeta|^4+\textbf{i}\langle\frac{\sbmx}{t^{1/4}},\sbmeta\rangle_N}\,d\bmeta. 
\end{equation}
Taking account of $\partial_{z_i}^2e^{\textbf{i}\langle\sbmz,\sbmeta\rangle_N}
=-\eta_i^2e^{\textbf{i}\langle\sbmz,\sbmeta\rangle_N}\,
(i=1,\cdots,N)$, we obtain the desired result by some properties of the Fourier transform and \eqref{rep_G_Nlambda}. \qed

\bigskip
Set $\mathbb{Z}_+^N:=\{\alpha=(\alpha_1,\cdots,\alpha_N)\in\mathbb{Z}^N\,|\,\alpha_i\ge0\,(i=1,\cdots,N)\}$ and 
$\mathbb{Z}_+:=\mathbb{Z}_+^1$. 
Hereafter, $\alpha=(\alpha_1,\cdots,\alpha_N)\in\mathbb{Z}_+^N$ is a multi-index with $|\alpha|=\alpha_1+\cdots+\alpha_N$ 
and $D_{\sbmx}^\alpha=\partial_{x_1}^{\alpha_1}\cdots\partial_{x_N}^{\alpha_N}$. 

For $G_{N,\lambda}$, we have the pointwise estimates as follows:

\begin{thm}\label{thm:Cui}
There exist $C,\,\nu,\,\mu,\,K>0$ such that for all $\alpha\in\mathbb{Z}_+^N$, $m\in\mathbb{Z}_+$, $\lambda\in\bR$ 
and $ (\bmx,t)\in\bR^N\times(0,\infty) $
\begin{align*}
&|D_{\sbmx}^\alpha(-\Delta_{\sbmx})^m G_{N,\lambda}(\bmx,t)| \\
&\le C\,\nu^{|\al|+2m}\,\Gamma\Bigl(\frac{|\alpha|+2m+N}4\Bigr)t^{-(N+|\alpha|+2m)/4}
\biggl(1+\frac{|\bmx|}{t^{1/4}}\biggr)^{(|\alpha|+2m)/3}e^{-\mu(|\sbmx|^4/t)^{1/3}+K|\lambda|^2t}.
\end{align*}
\end{thm}

\noindent
In the case $ \lambda=0 $ this estimate is originally obtained in Eidel'man\,\cite[Section 3 in Chapter I]{eid;69} 
and Cui\,\cite[Theorem 3.2]{cui;01}. 
However, in these references, the precise dependence of the constant on $\alpha\in\mathbb{Z}_+^N$ and $m\in\mathbb{Z}_+$ 
is not stated. Thus we give the proof of this theorem in Appendix B below. 
The following corollary is a direct consequence of Theorem \ref{thm:Cui}. 

\begin{cor}\label{cor:Cui-thm}
There exist $C,\,\nu,\,\mu>0$ such that for all $\alpha\in\mathbb{Z}_+^N$, $m\in\mathbb{Z}_+$ and $ \bmx\in\bR^N $
\[
|D_{\sbmx}^\alpha(-\Delta_{\sbmx})^m g_N(\bmx)|
\le C\,\nu^{|\al|+2m}\,\Gamma\Bigl(\frac{|\alpha|+2m+N}4\Bigr)(1+|\bmx|)^{(|\alpha|+2m)/3}
e^{-\mu|\sbmx|^{4/3}}.
\]
\end{cor}

\medskip\noindent
If $N\ge2$, for an orthogonal matrix $P_N$ of size $N$ satisfying $P_N\bmx=(|\bmx|,0,\cdots,0)^T\in\mathbb{R}^N$, 
set $\bmxi=P_N\bmzeta$ in \eqref{g_N}. Since $|\det P_N\,|=1$ and $|\bmxi|=|\bmzeta|$, we obtain
\[
g_N(\bmx)=c_N\int_{\mathbb{R}^N}e^{-|\sbmzeta|^4+\textbf{i}\langle\sbmx,P_N\sbmzeta\rangle_N}\,d\bmzeta
=c_N\int_{\mathbb{R}^N}e^{-|\sbmzeta|^4+\textbf{i}|\sbmx|\zeta_1}\,d\bmzeta.
\]
Applying the change of variable on the polar coordinate, $g_N(\bmx)$ is represented as 
\begin{align}\label{rot_sym1}
g_N(\bmx)
=&\,c_N\omega_{N-2}\int_0^{\infty}\rho^{N-1}e^{-\rho^4}
\left\{\int_0^{\pi}e^{\textbf{i}|\sbmx|\rho\cos\theta}\sin^{N-2}\ms\theta\,d\theta\right\}d\rho \\
=&\,(2\pi)^{-N/2}|\bmx|^{1-N}\int_0^{\infty}(|\bmx|\rho)^{N/2}e^{-\rho^4}J_{(N-2)/2}(|\bmx|\rho)\,d\rho \nonumber
\end{align}
where  $J_{(N-2)/2}$ is the $(N-2)/2$-th Bessel function. If $N=1$, we have
\begin{align}\label{rot_sym2}
g_1(x)
=&\,c_1\int_{\mathbb{R}}e^{-\xi^4+\textbf{i}x\xi}\,d\xi=2c_1\int_0^\infty e^{-\xi^4}\cos(x\xi)\,d\xi \\
=&\,(2\pi)^{-1/2}\int_0^{\infty}(|x|\xi)^{1/2}e^{-\xi^4}J_{-1/2}(|x|\xi)\,d\xi. \nonumber
\end{align}
From \eqref{rot_sym1} and \eqref{rot_sym2}, set
\[
\phi_N(|\bmx|):=(2\pi)^{-N/2}|\bmx|^{1-N}\int_0^{\infty}(|\bmx|\rho)^{N/2}e^{-\rho^4}J_{(N-2)/2}(|\bmx|\rho)\,d\rho
\]
for $N\ge1$, which means $g_N(\bmx)=\phi_N(|\bmx|)$. 
According to Ferrero--Gazzola--Grunau\,\cite[Section 2]{fer;gaz;gru;08}, $\phi_N(r)$ with $r=|\bmx|\ge0$ satisfies 
\begin{align}
&\phi_N(r)
=\frac1{2^{N+1}\pi^{N/2}}\sum_{\ell=0}^\infty
\frac{(-1)^\ell\Gamma(\ell/2+N/4)}{2^{2\ell}\Gamma(\ell+1)\Gamma(\ell+N/2)}r^{2\ell}, 
\label{phi_N} \\
&\phi_N'(r)=-r\phi_{N+2}(r). \label{1ode-phi_N}
\end{align}
In addition, we have the following lemma. 

\begin{lem}\label{lem:m-derivative_gN}
For $g_N(\bmx)=\phi_N(r)$ with $r=|\bmx|\ge0$,  
\begin{align*}
\Delta^mg_N(\bmx)
=&\,\left(\frac{d^2}{dr^2}+\frac{N-1}{r}\frac{d}{dr}\right)^m\phi_N(r) \\
=&\,\frac{(-1)^m}{2^{N+1}\pi^{N/2}}\sum_{\ell=0}^\infty
\frac{(-1)^\ell\Gamma((\ell+m)/2+N/4)}{2^{2\ell}\Gamma(\ell+1)\Gamma(\ell+N/2)}\,r^{2\ell}.
\end{align*}
\end{lem}
\noindent
{\bf Proof.} 
For $m=1$, the direct calculation yields that
\begin{align*}
\left(\frac{d^2}{dr^2}+\frac{N-1}{r}\frac{d}{dr}\right)\phi_N(r)
& =\frac{-1}{2^{N+1}\pi^{N/2}}\sum_{\ell=0}^\infty 
\frac{(-1)^\ell\Gamma((\ell+1)/2+N/4)}{2^{2\ell}\Gamma(\ell+1)\Gamma(\ell+N/2)}\,r^{2\ell}.
\end{align*}
The result follows by induction. \qed

\bigskip

We show the following lemma.  It is necessary to estimate some integrations of $ G_{N,\lambda} $ 
in the next section. 

\begin{lem}\label{lem:series_finite}
Let $\alpha\in\mathbb{Z}_+^N$ with $|\alpha|\le2$. Then there exists $\gamma>0$ such that for $|\lambda|t^{1/2}\le\gamma$
\[
\sum_{m=0}^\infty\frac{|\lambda|^mt^{m/2}}{m!}\int_{\mathbb{R}^N}\bigl|
D_{\sbmx}^\alpha(-\Delta_{\sbmx})^mg_N(\bmx)\bigr|\,d\bmx<\infty.
\]
\end{lem}
\noindent
{\bf Proof.} 
By Corollary \ref{cor:Cui-thm}, there exist $C,\,\nu,\,\mu>0$ independent of $m$, such that 
\begin{align*}
&\frac{|\lambda|^m t^{m/2}}{m!}\int_{\mathbb{R}^N}
\bigl|D_{\sbmx}^\alpha(-\Delta_{\sbmx})^mg_N(\bmx)\bigr|\,d\bmx \\
&\le\frac{C\,(|\lambda|\nu^2t^{1/2})^m}{m!}
\Gamma\Bigl(\frac{m}2+\frac{|\alpha|+N}4\Bigr)\int_{\mathbb{R}^N}
(1+|\bmx|)^{(|\alpha|+2m)/3}e^{-\mu|\sbmx|^{4/3}}\,d\bmx.
\end{align*}
Applying the change of variable on the polar coordinate, we obtain 
\begin{align*}
&\hspace*{-12pt}
\int_{\mathbb{R}^N}(1+|\bmx|)^{(|\alpha|+2m)/3}e^{-\mu|\sbmx|^{4/3}}\,d\bmx \\
\le&\,2^{(|\alpha|+2m)/3}C_1\left\{1+\left(\frac1{\mu}\right)^{(|\alpha|+2m)/4}
\Gamma\Bigl(\frac{m}2+\frac{|\alpha|+3N}4\Bigr)\right\}
\end{align*}
for a constant $C_1>0$ independent of $m$. Here it follows from the Schwarz inequality that 
\[
\Gamma\Bigl(\frac{m}2+\frac{|\alpha|+kN}4\Bigr)
\le\{\Gamma(m)\}^{1/2}\biggl\{\Gamma\Bigl(\frac{|\alpha|+kN}2\Bigr)\biggr\}^{1/2}
\le(m!)^{1/2}\biggl\{\Gamma\Bigl(\frac{|\alpha|+kN}2\Bigr)\biggr\}^{1/2}
\]
for $m\in\mathbb{N}$ and $k=1,3$. Therefore, we see that 
\[
\frac{|\lambda|^mt^{m/2}}{m!}\int_{\mathbb{R}^N}
\bigl|D_{\sbmx}^\alpha(-\Delta_{\sbmx})^mg_N(\bmx)\bigr|\,d\bmx
\le C_2\,\frac{(4^{1/3}|\lambda|\nu^2t^{1/2})^m}{(m!)^{1/2}}
+C_3\left(\frac{4^{1/3}|\lambda|\nu^2t^{1/2}}{\mu^{1/2}}\right)^m
\]
for constants $C_2,C_3>0$ independent of $m$. Choose $ \gamma>0 $ satisfying 
\[
\gamma<\frac{\mu^{1/2}}{\,4^{1/3}\nu^2\,}.
\]
Then, with the help of d'Alembert test, we can judge that the series
\[
\sum_{m=1}^\infty\frac{(4^{1/3}|\lambda|\nu^2t^{1/2})^m}{(m!)^{1/2}}, \quad 
\sum_{m=1}^\infty\left(\frac{4^{1/3}|\lambda|\nu^2t^{1/2}}{\mu^{1/2}}\right)^m
\]
converge uniformly for $|\lambda|t^{1/2}\le\gamma$, so that we obtain the desired result. \qed

\bigskip\noindent
Lemma \ref{lem:series_finite} and the Lebesgue convergence theorem lead to the following lemma. 

\begin{lem}\label{lem:termwise-int}
Let $\alpha\in\mathbb{Z}_+^N$ with $|\alpha|\le2$ and $h\in L^\infty(\mathbb{R}^N)$. 
Then there exists $\gamma>0$ such that for $|\lambda|t^{1/2}\le\gamma$
\begin{align*}
&\int_{\mathbb{R}^N}\sum_{m=0}^\infty\frac{(-\lambda)^mt^{m/2}}{m!}D_{\sbmx}^\alpha(-\Delta_{\sbmx})^mg_N(\bmx)\,h(\bmx)\,d\bmx \\
&=\sum_{m=0}^\infty\frac{(-\lambda)^mt^{m/2}}{m!}\int_{\mathbb{R}^N}D_{\sbmx}^\alpha(-\Delta_{\sbmx})^mg_N(\bmx)\,h(\bmx)\,d\bmx.  
\end{align*}
\end{lem}
 
\medskip\noindent
With regard to $(-\partial_{x_N}^2)^\ell(-\Delta_{\sbmx})^mg_N(\bmx',0)$, we have the following representation.
 
\begin{lem}\label{lem:g_xN=0}
Let $g_N$ be given by \eqref{g_N}. Then 
\begin{align*}
&(-\partial_{x_N}^2)^\ell(-\Delta_{\sbmx})^mg_N(\bmx',0) \\
&=c_1\sum_{j=0}^m\binom{m}{j}\mathcal{F}^{-1}_{N-1}
\biggl[\,\sum_{k=0}^\infty\frac{(-2)^k}{k!}L_{2(k+j+\ell)}|\,\cdot\,|^{2(k+m-j)}e^{-|\,\cdot\,|^4}\,\biggr](\bmx')
\end{align*}
for $\ell,m\in\mathbb{Z}_+$, where
\begin{equation}\label{cofficient-L}
L_\sigma:=2\int_0^\infty\xi^\sigma e^{-\xi^4}\,d\xi=\frac12\Gamma\Bigl(\frac{\sigma+1}4\Bigr)
\end{equation}
for $\sigma\ge0$. 
\end{lem}
\noindent
{\bf Proof.} 
Fix any $\bmx\in\mathbb{R}^N$. 
Taking account of $|\bmxi|^4=|\bmxi'|^4+\xi_N^4+2|\bmxi'|^2\xi_N^2$ for $\bmxi'=(\xi_1,\cdots,\xi_{N-1})$, we see that 
for $\bmx'=(x_1,\cdots,x_{N-1})$
\begin{align*}
&\hspace*{-10pt}
(-\partial_{x_N}^2)^\ell(-\Delta_{\sbmx})^mg_N(\bmx) \\
=&\,c_N\int_{\mathbb{R}^N}\xi_N^{2\ell}|\bmxi|^{2m}e^{-|\sbmxi|^4+\textbf{i}\langle\sbmx,\sbmxi\rangle_N}\,
d\bmxi \nonumber \\
=&\,c_N\sum_{j=0}^m\binom{m}{j}\int_{\mathbb{R}^N}|\bmxi'|^{2(m-j)}\xi_N^{2(j+\ell)}
e^{-(|\sbmxi'|^4+\xi_N^4+2|\sbmxi'|^2\xi_N^2)+\textbf{i}\{\langle\sbmx',\sbmxi'\rangle_{N-1}+x_N\xi_N\}}\,
d\bmxi \nonumber \\
=&\,c_N\sum_{j=0}^m\binom{m}{j}\int_{\mathbb{R}^N}\sum_{k=0}^{\infty}\frac{\,(-2)^k\,}{k!}|
\bmxi'|^{2(k+m-j)}\xi_N^{2(k+j+\ell)}
e^{-|\sbmxi'|^4+\textbf{i}\langle\sbmx',\sbmxi'\rangle_{N-1}}e^{-\xi_N^4+\textbf{i}x_N\xi_N}\,d\bmxi.
\nonumber
\end{align*}
Note that by Fubini's theorem
\begin{align*}
(-\partial_{x_N}^2)^\ell(-\Delta_{\sbmx})^mg_N(\bmx)
=&\,c_N\sum_{j=0}^m\binom{m}{j}
\int_{\mathbb{R}^{N-1}}|\bmxi'|^{2(m-j)}e^{-|\sbmxi'|^4+\textbf{i}\langle\sbmx',\sbmxi'\rangle_{N-1}} \\
&\hspace*{50pt}
\left\{\int_{\mathbb{R}}\sum_{k=0}^\infty\frac{(-2)^k}{k!}
|\bmxi'|^{2k}\xi_N^{2(k+j+\ell)}e^{-\xi_N^4+\textbf{i}\,x_N\xi_N}\,d\xi_N\right\}d\bmxi'.
\end{align*}
Here it follows that for each $\bmxi'\in\mathbb{R}^{N-1}$ and $n\in\mathbb{N}$, 
\begin{align*}
&\int_{\mathbb{R}}\left|\sum_{k=0}^n\frac{(-2)^k}{k!}|\bmxi'|^{2k}\xi_N^{2(k+j+\ell)}e^{-\xi_N^4+\textbf{i}\,x_N\xi_N}\right|d\xi_N
\le\int_{\mathbb{R}}\xi_N^{2(j+\ell)}e^{-\xi_N^4+2|\sbmxi'|^2\xi_N^2}d\xi_N<+\infty.
\end{align*}
The Lebesgue convergence theorem implies that
\[
\int_{\mathbb{R}}\sum_{k=0}^\infty\frac{(-2)^k}{k!}|\bmxi'|^{2k}\xi_N^{2(k+j+\ell)}e^{-\xi_N^4+\textbf{i}\,x_N\xi_N}\,d\xi_N
=\sum_{k=0}^\infty\frac{(-2)^k}{k!}|\bmxi'|^{2k}\int_{\mathbb{R}}\xi_N^{2(k+j+\ell)}e^{-\xi_N^4+\textbf{i}\,x_N\xi_N}\,d\xi_N.
\]
Substituting $x_N=0$, we have
\begin{align*}
\int_{\mathbb{R}}\sum_{k=0}^\infty\frac{(-2)^k}{k!}|\bmxi'|^{2k}\xi_N^{2(k+j+\ell)}e^{-\xi_N^4}\,d\xi_N
=&\,\sum_{k=0}^\infty\frac{(-2)^k}{k!}|\bmxi'|^{2k}\int_{\mathbb{R}}\xi_N^{2(k+j+\ell)}e^{-\xi_N^4}\,d\xi_N \\
=&\,\sum_{k=0}^\infty\frac{(-2)^k}{k!}|\bmxi'|^{2k}L_{2(k+j+\ell)}.
\end{align*}
Since 
\[
|\,\cdot\,|^{2(m-j)}e^{-|\,\cdot\,|^4+\textbf{i}\langle\sbmx',\,\cdot\,\rangle_{N-1}}\sum_{k=0}^\infty
\frac{(-2)^k}{k!}L_{2(k+j+\ell)}|\,\cdot\,|^{2k}\in L^1(\mathbb{R}^{N-1})
\]
by Fubini's Theorem, we obtain the desired result. \qed

\section{Asymptotic expansion of a solution to linear parabolic equations} \label{sec:asympt-exp_4w2}
The purpose of this section is to derive the asymptotic expansion of a solution $u(\bmx,t)$ to \eqref{4th_diffusion} as 
$t\to+0$. 
Throughout this section, we assume that $\Omega_0$ is a compact set in $\mathbb{R}^N$ and $\partial \Omega_0$ 
is of class $C^5$. 
Recalling Proposition \ref{prop:rep_4th-laplace} and Lemma \ref{lem:termwise-int}, the solution $u(\bmx,t)$ 
to \eqref{4th_diffusion} is given by 
\begin{align*}
u(\bmx,t)
=&\,(G_{N,\lambda}(\cdot,t)\ast\chi_{\Omega_0})(\bmx) \\
=&\,\frac1{t^{N/4}}\sum_{m=0}^\infty\frac{(-\lambda)^mt^{m/2}}{m!}
\int_{\mathbb{R}^N}(-\Delta_{\sbmz})^mg_N\Bigl(\frac{\bmx-\bmy}{t^{1/4}}\Bigr)\chi_{\Omega_0}(\bmy)\,d\bmy \\
=&\,\sum_{m=0}^\infty\frac{(-\lambda)^mt^{m/2}}{m!}
\int_{\{\sbmz\in\mathbb{R}^N\,|\,\sbmx-t^{1/4}\sbmz\in \Omega_0\}}
(-\Delta_{\sbmz})^mg_N(\bmz)\,d\bmz.
\end{align*}
for $\bmx\in\mathbb{R}^N$. 

\subsection{Notations related to $\p\om_0$} \label{sec:geometric-asmption}
In this subsection, we give some notations related to $\p\om_0$. 
For $\bmx_\ast\in\mathbb{R}^N$ and $\delta>0$, we denote a neighborhood of $\bmx_\ast$ by 
\begin{align*}
&Q_{\sbmx_\ast,\delta}:=\{\bmx\in\mathbb{R}^N\,|\,|x_i-x_{\ast,i}|<\delta\,(i=1,\cdots,N)\}, \quad 
Q_{\sbmx_\ast}:=Q_{\sbmx_\ast,1}, \\
&Q'_{\sbmx_\ast,\delta}:=\{\bmx\in Q_{\sbmx_\ast,\delta}\,|\,x_N=x_{\ast,N}\}, \quad 
Q'_{\sbmx_\ast}:=Q'_{\sbmx_\ast,1}. 
\end{align*}
Since $ \partial\Omega_0 $ is compact, we can choose $ \delta_0>0 $ such that the family 
$ \{Q_{\sbmx_\ast,\delta_0}\}_{\sbmx_\ast\in\partial\Omega_0} $ is an open covering of 
$ \partial\Omega_0 $. We may set $ \delta_0=1 $. 
Let $\bmn(\bmx)$ be the unit outward normal to $\partial\Omega_0$ at $\bmx\in\partial\Omega_0$ and 
$P_N$ be an orthogonal matrix of size $N$ such that $P_N\bmn(\bmx)=\bme_N$ where 
$\bme_N=(\bmzero',1)^T$ with $\bmzero'=(0,\cdots,0)\in\mathbb{R}^{N-1}$. Then, set 
$\widetilde{\Omega}_0:=\{\widetilde{\bmy}\in\mathbb{R}^N\,|\,P_N^{-1}\widetilde{\bmy}\in \Omega_0\}$. 
Without loss of generality, we may assume that for each $\bmx\in\partial\Omega_0$ there exists 
a function $f|_{\sbmx}:Q'_{\sbmzero}\to\mathbb{R}$ satisfying the following properties:
\begin{list}{}{\topsep=0.15cm\itemsep=0cm\leftmargin=1.2cm\labelwidth=1cm}
\item[(A1)]
$f|_{\sbmx}\in C^5(\overline{Q'_{\sbmzero}})$ for any $\bmx\in\p\om_0$ and 
$\|f|_{\sbmx}\|_{C^5(\overline{Q'_{\sbmzero}})}$ is uniformly bounded for $\bmx\in\p\om_0$.
\item[(A2)]
$f|_{\sbmx}(\bmzero')=0$ and $\nabla_{\sbmx^\prime} f|_{\sbmx}(\bmzero')=\bmzero'$. 
\item[(A3)]
$\partial\widetilde{\Omega}_0\cap Q_{P_N\sbmx}=\{\widetilde{\bmy}\in\mathbb{R}^N\,|\,
\widetilde{y}_N-(P_N\bmx)_N=f|_{\sbmx}(\widetilde{\bmy}'-(P_N\bmx)')\,\,(\widetilde{\bmy}'\in Q'_{P_N\sbmx})\}$ 
where $(P_N\bmx)_i$ is the $i$-th component of $P_N\bmx$ and $(P_N\bmx)'=((P_N\bmx)_1,\cdots,(P_N\bmx)_{N-1})$. 
\end{list}
Note that $(\widetilde{\bmy}',(P_N\bmx)_N+f|_{\sbmx}(\widetilde{\bmy}'-(P_N\bmx)'))\,\,(\widetilde{\bmy}'\in Q'_{P_N\sbmx})$ 
is a graph representation of $\partial\widetilde{\Omega}_0$ in a neighborhood of  $P_N\bmx$ 
for $\bmx\in\partial\Omega_0$. Hereafter, for simplicity, we denote $f|_{\sbmx}$ by $f$. 
We also define a function $\psi=\psi(\bmz',v,t)\colon\bR^{N-1}\times\bR\times \bR_{+}\to\bR$ as
\[
\psi(\bmz',v,t):=t^{-1/4}\{-v+f(t^{1/4}\bmz')\}.
\]
Let $g=(g_{ij})$, $A=(h_{ij})$, $H$, $\kappa_i$ and $\Delta_g$ be the induced metric, the second fundamental form, 
the mean curvature, 
the principal curvatures and the Laplace-Beltrami operator of $\partial\Omega_0$, respectively.
$|A|^2$ is the norm of the second fundamental form, that is, $|A|^2=\sum\limits_{i,j,k,\ell\in\Lambda}g^{ij}g^{k\ell}h_{ik}h_{j\ell}$, 
where $g^{-1}=(g^{ij})$ is the inverse matrix of $g=(g_{ij})$. 
For the representation of these quantities by $f$, see Appendix \ref{sec:mc}. 

Define the signed distance function to $\partial\Omega_0$ as 
\begin{equation}\label{signed-distance}
d(\bmy,\partial \Omega_0):=\left\{\begin{array}{rl} 
\inf\limits_{\sbmx\in\partial \Omega_0}|\bmy-\bmx|&(\bmy\in\Omega_0), \\ 
-\inf\limits_{\sbmx\in\partial \Omega_0}|\bmy-\bmx|&(\bmy\in\mathbb{R}^N\setminus\Omega_0). 
\end{array}\right.
\end{equation}
Set $(\partial\Omega_0)^\delta:=\{\bmy\in\mathbb{R}^N\,|\,|d(\bmy,\partial\Omega_0)|<\delta\}$. 
Then we take $\delta_0\in(0,1/2)$ such that for any $\bmy\in(\partial\Omega_0)^{\delta_0}$ there is a unique 
$\bmx_0\in\partial\Omega_0$ satisfying $|d(\bmy,\partial\Omega_0)|=|\bmy-\bmx_0|$. 

Furthermore, in the following, we often use the notation $F(\bmx,t)=O_D(b(t))$. This means that there exists 
$C>0$, which is obtained uniformly for $\bmx\in D$, such that $|F(\bmx,t)|\le C|b(t)|$. 

\subsection{Asymptotic expansion}
We first show the asymptotic expansion of $u(\bmx+v\,\bmn(\bmx),t)$ by using the graph representation: 

\begin{thm}\label{thm:asympt-exp_4w2}
Let $\gamma$ be a constant obtained in Lemma \ref{lem:series_finite}. Then, there exists $\mu_*>0$ such that for 
all $\bmx\in\partial \Omega_0$, $v\in(-\delta_0,\delta_0)$, $\lambda\in\mathbb{R}$ and $t>0$ satisfying 
$|\lambda|t^{1/2}\le\gamma$ 
\begin{align}\label{sol_expansion}
&\hspace*{-10pt}
u(\bmx+v\,\bmn(\bmx),t) \\
=&\,\frac12+\sum_{m=0}^\infty\frac{(-\lambda)^mt^{m/2}}{m!}
\int_{\mathbb{R}^{N-1}}\int_0^{\psi(\sbmz',v,t)}(-\Delta_{\sbmz})^mg_N(\bmz)\,dz_Nd\bmz'
+O_{\p\om_0}(e^{-\mu_\ast t^{-1/3}}). \nonumber
\end{align}

Furthermore, for $j\in\{1,\cdots, N\}$, $\lambda\in\mathbb{R}$ and $t>0$ satisfying
$|\lambda|t^{1/2}\le\gamma$
\begin{align}\label{Dsol_expansion}
&\hspace*{-10pt}
\langle\nabla_{\sbmx} u(\bmx+v\,\bmn(\bmx),t), P_N^{-1}\bme_j\rangle_N \\
=&\,-\frac1{t^{1/4}}\sum_{m=0}^\infty\frac{(-\lambda)^mt^{m/2}}{m!}
\int_{\mathbb{R}^{N-1}}\int_{-\infty}^{\psi(\sbmz',v,t)}\partial_{z_j}(-\Delta_{\sbmz})^mg_N(\bmz)\,dz_Nd\bmz' 
\nonumber \\
&\,+O_{\p\om_0}(e^{-\mu_\ast t^{-1/3}}). \nonumber
\end{align}

\end{thm}

\begin{rem}
Note that the term $ O_{\p\om_0}(e^{-\mu_\ast t^{-1/3}}) $ is not only uniform 
for $ \bmx\in\partial\Omega_0 $, but also for $ v\in(-\delta_0,\delta_0) $.  See the estimate of $ I_2 $ in the proof below. 
\end{rem}

\medskip\noindent
{\bf Proof of Theorem \ref{thm:asympt-exp_4w2}.} 
{\it Step 1.} We first prove (\ref{sol_expansion}). 
Using the fact that $G_{N,\lambda}(\,\cdot\,,t)$ is radially symmetric, we have
\begin{align*}
u(\bmx+v\,\bmn(\bmx),t)
=&\,\int_{\Omega_0}G_{N,\lambda}(\bmx+v\,\bmn(\bmx)-\bmy,t)\,d\bmy \\
=&\,\int_{\Omega_0}G_{N,\lambda}(\bmy-(\bmx+v\,\bmn(\bmx)),t)\,d\bmy.
\end{align*}
Choose an orthogonal matrix $P_N$ of size $N$ such that $P_N\bmn(\bmx)=\bme_N$. Setting 
$\widetilde{\bmy}=P_N\bmy$ 
and taking account of $|P_N^{-1}\bmz|=|\bmz|$ for $\bmz\in\mathbb{R}^N$, we obtain
\begin{align*}
u(\bmx+v\,\bmn(\bmx),t)
&=\int_{\widetilde{\Omega}_0}G_{N,\lambda}(\widetilde{\bmy}-(P_N\bmx+v\,\bme_N),t)\,d\widetilde{\bmy} \\
&=\left(\int_{\widetilde{\Omega}_0\cap Q_{{\tiny P_ N}\tbmx}}
+\int_{\widetilde{\Omega}_0\setminus Q_{{\tiny P_N}\tbmx}}\right)
G_{N,\lambda}(\widetilde{\bmy}-(P_N\bmx+v\,\bme_N),t)\,d\widetilde{\bmy} \\
&=:I_1+I_2.
\end{align*}
It follows from Theorem \ref{thm:Cui} with $ |\al|=m=0 $ that for $t>0$ satisfying $|\lambda|t^{1/2}\le\gamma$ 
\begin{align*}
I_2\le&\,C_1\,t^{-N/4}e^{K|\lambda|^2t}\int_{\mathbb{R}^N\setminus B_N(P_N\sbmx+v\,\sbme_N,1/2)}
e^{-\mu(|\tilde{\sbmy}-(P_N\sbmx+v\,\sbme_N)|^4/t)^{1/3}}\,d\widetilde{\bmy} \\
\le&\,C_2\int_{\mathbb{R}^N\setminus B_N(\sbmzero,t^{-1/4}/2)}e^{-\mu|\tilde{\sbmeta}|^{4/3}}\,d\widetilde{\bmeta}
=O_{\p\om_0}(e^{-\mu_\ast t^{-1/3}}),
\end{align*}
where $B_N(\bmx,r)$ is the $N$-dimensional ball with the center $\bmx$ and the radius $r$. We also see that for 
$t>0$ satisfying $|\lambda|t^{1/2}\le\gamma$ 
\begin{align*}
&\left|\int_{\mathbb{R}^{N-1}}\int_{-\infty}^{f(\tilde{\sbmy}'-(P_N\sbmx)')+(P_N\sbmx)_N}
G_{N,\lambda}(\widetilde{\bmy}-(P_N\bmx+v\,\bme_N),t)\,d\tilde{y}_Nd\widetilde{\bmy}'-I_1\right| \\
&\le C\,t^{-N/4}e^{K|\lambda|^2t}\int_{\mathbb{R}^N\setminus B_N(P_N\sbmx+v\,\sbme_N,1/2)}
e^{-\mu(|\tilde{\sbmy}-(P_N\sbmx+v\,\sbme_N)|^4/t)^{1/3}}\,d\widetilde{\bmy}
=O_{\p\om_0}(e^{-\mu_\ast t^{-1/3}}),
\end{align*}
where $f$ is a function satisfying (A1)-(A3). These facts imply that 
\begin{align*}
u(\bmx+v\,\bmn(\bmx),t)
=&\,\int_{\mathbb{R}^{N-1}}\int_{-\infty}^{f(\tilde{\sbmy}'-(P_N\sbmx)')+(P_N\sbmx)_N}
G_{N,\lambda}(\widetilde{\bmy}-(P_N\bmx+v\,\bme_N),t)\,d\tilde{y}_Nd\widetilde{\bmy}' \\
&\,+O_{\p\om_0}(e^{-\mu_\ast t^{-1/3}}).
\end{align*}
Applying the change of variable $\bmeta=\widetilde{\bmy}-(P_N\bmx+v\,\bme_N)$ and recalling \eqref{integralG=1}, 
Proposition \ref{prop:rep_4th-laplace} and Lemma \ref{lem:termwise-int}, we have
\begin{align*}
&\hspace*{-10pt}
u(\bmx+v\,\bmn(\bmx),t) \\
=&\,\int_{\mathbb{R}^{N-1}}\int_{-\infty}^{\psi(\sbmeta',v,1)}G_{N,\lambda}(\bmeta,t)\,d\eta_Nd\bmeta'
+O_{\p\om_0}(e^{-\mu_\ast t^{-1/3}}) \\
=&\,\frac12+\int_{\mathbb{R}^{N-1}}\int_0^{\psi(\sbmeta',v,1)}G_{N,\lambda}(\bmeta,t)\,d\eta_Nd\bmeta'
+O_{\p\om_0}(e^{-\mu_\ast t^{-1/3}}) \\
=&\,\frac12+\frac1{t^{N/4}}\sum_{m=0}^\infty\frac{(-\lambda)^mt^{m/2}}{m!}
\int_{\mathbb{R}^{N-1}}\int_0^{\psi(\sbmeta',v,1)}(-\Delta_{\sbmz})^m
 g_N\Bigl(\frac{\bmeta}{t^{1/4}}\Bigr)\,d\eta_Nd\bmeta' \\
&\,+O_{\p\om_0}(e^{-\mu_\ast t^{-1/3}}).
\end{align*}
By the change of variable $\bmz=t^{-1/4}\bmeta$, we obtain \eqref{sol_expansion}. 

{\it Step 2.} We derive \eqref{Dsol_expansion}.  
Using again the fact that $G_{N,\lambda}(\,\cdot\,,t)$ is radially symmetric and setting $\widetilde{\bmy}=P_N\bmy$ 
with an orthogonal matrix $P_N$ introduced in Step 1, we see that 
\begin{align*}
\nabla_{\sbmx} u(\bmx+v\,\bmn(\bmx),t)
=&\,-\int_{\widetilde{\Omega}_0}P_N^T\nabla_{\sbmz}G_{N,\lambda}(\widetilde{\bmy}-(P_N\bmx+v\,\bme_N),t)\,d\widetilde{\bmy},
\end{align*}
where $P_N^T$ is the transposed matrix of $P_N$ and 
$\nabla_{\sbmz}G_{N,\lambda}=(\partial_{z_1}G_{N,\lambda}(\bmz,t),\cdots,\partial_{z_N}G_{N,\lambda}(\bmz,t))$. 
This implies that
\begin{align*}
\langle\nabla_{\sbmx} u(\bmx+v\,\bmn(\bmx),t),P_N^{-1}\bme_{j}\rangle_N
&=\,-\int_{\widetilde{\Omega}_0}\langle P_N^T\nabla_{\sbmz}G_{N,\lambda}(\widetilde{\bmy}-(P_N\bmx+v\,\bme_N),t),P_N^{-1}\bme_{j}\rangle_N
\,d\widetilde{\bmy}\\
&=\,-\int_{\widetilde{\Omega}_0}\partial_{z_j}G_{N,\lambda}(\widetilde{\bmy}-(P_N\bmx+v\,\bme_N),t)\,d\widetilde{\bmy}.
\end{align*}
Applying an argument similar to the above, we are led to \eqref{Dsol_expansion}. \qed

\bigskip
By the Taylor's expansion of $f(\bmz')$ at $\bmz'=\bmzero'$, we see that 
\begin{align}\label{surface_expansion}
\psi(\bmz',v,t)
=&\,t^{-1/4}\biggl\{
-v+f(\bmzero')+\langle\bmz',\nabla_{\sbmz'}\rangle_{N-1}f(\bmzero')t^{1/4}
+\frac{\langle\bmz',\nabla_{\sbmz'}\rangle_{N-1}^2f(\bmzero')}2t^{1/2} \\
&\hspace*{35pt}
+\frac{\langle\bmz',\nabla_{\sbmz'}\rangle_{N-1}^3f(\bmzero')}6t^{3/4}
+\frac{\langle\bmz',\nabla_{\sbmz'}\rangle_{N-1}^4f(\bmzero')}{24}t \nonumber \\
&\hspace*{35pt}
+\frac{\langle\bmz',\nabla_{\sbmz'}\rangle_{N-1}^5f(\theta t^{1/4}\bmz')}{120}t^{5/4}\biggr\} \nonumber \\
=&\,-vt^{-1/4}+\frac{\langle\bmz',\nabla_{\sbmz'}\rangle_{N-1}^2f(\bmzero')}2t^{1/4}
+\frac{\langle\bmz',\nabla_{\sbmz'}\rangle_{N-1}^3f(\bmzero')}6t^{1/2} \nonumber \\
&\,+\frac{\langle\bmz',\nabla_{\sbmz'}\rangle_{N-1}^4f(\bmzero')}{24}t^{3/4}
+\frac{\langle\bmz',\nabla_{\sbmz'}\rangle_{N-1}^5f(\theta t^{1/4}\bmz')}{120}t \nonumber
\end{align}
for some $\theta\in(0,1)$. Note that $ \langle\bmz',\nabla_{\sbmz'}\rangle_{N-1}^nf $ is 
defined by 
\begin{align*}
\langle\bmz',\nabla_{\sbmz'}\rangle_{N-1}^nf
=&\,\left(\sum_{i\in\Lambda}z_i\partial_{z_i}\right)^nf \\
=&\,\sum_{d_1+\cdots+d_{N-1}
=n}\binom{n}{d_1,\cdots,d_{N-1}}(z_1\partial_{z_1})^{d_1}\cdots(z_{N-1}\partial_{z_{N-1}})^{d_{N-1}}f
\end{align*}
where 
\[
\binom{n}{d_1,\cdots,d_{N-1}}=\frac{n!}{d_1!\cdots d_{N-1}!}.
\]
According to (A2) and Appendix \ref{sec:mc} below, $\langle\bmz',\nabla_{\sbmz'}\rangle_{N-1}^2f(\bmzero')$ 
is represented as 
\begin{equation}\label{kappa1}
\langle\bmz',\nabla_{\sbmz'}\rangle_{N-1}^2f(\bmzero')=\sum_{i\in\Lambda}\kappa_i\zeta_i^2,
\end{equation}
where 
$\bmzeta'=P_{N-1}^{-1}\bmz'$ for an orthogonal matrix $P_{N-1}$ of size $N-1$. 
This implies 
\begin{align}
\bigl(\langle\bmz',\nabla_{\sbmz'}\rangle_{N-1}^2f\bigr)^2(\bmzero')
=&\,\sum_{i\in\Lambda}\kappa_i^2\zeta_i^4
+2\sum_{\substack{i_1,i_2\in\Lambda \\ i_1<i_2}}\kappa_{i_1}\kappa_{i_2}\zeta_{i_1}^2\zeta_{i_2}^2, \label{kappa2} \\
\bigl(\langle\bmz',\nabla_{\sbmz'}\rangle_{N-1}^2f\bigr)^3(\bmzero')
=&\,\sum_{i\in\Lambda}\kappa_i^3\zeta_i^6+3\sum_{\substack{i_1,i_2\in\Lambda \\ i_1\ne i_2}}
\kappa_{i_1}^2\kappa_{i_2}\zeta_{i_1}^4\zeta_{i_2}^2
+6\sum_{\substack{i_1,i_2,i_3\in\Lambda \\ i_1<i_2<i_3}}
\kappa_{i_1}\kappa_{i_2}\kappa_{i_3}\zeta_{i_1}^2\zeta_{i_2}^2\zeta_{i_3}^2. 
\label{kappa3}
\end{align}
Since $P_{N-1}$ is an orthogonal matrix, for any $F\colon\bR^{N-1}\to\bR$ it holds that
\[
	\int_{\mathbb{R}^{N-1}}F(\bmz')\,d\bmz'=\int_{\mathbb{R}^{N-1}}\widetilde{F}(\bmzeta')\,d\bmzeta',
\]
where $\widetilde{F}(\bmzeta'):=F(P_{N-1}\bmzeta')=F(\bmz')$. 
Hereafter, we use the variable $\bmz'$ in the sense of  the variable $\bmzeta'$ as the above. 

We introduce some notations. Set 
\[
B:=\{i\in\Lambda\,|\,\beta_i\ne0\}.
\] 
Then, for $\beta=(\beta_1,\cdots,\beta_{N-1})\in\mathbb{Z}_+^{N-1}$, $(\beta_i)_{i\in B}$ denotes a multi-index 
in which only positive components are chosen. For example, if $\beta_{i_1}=p\,(\,>0)$, $\beta_{i_2}=q\,(\,>0)$ 
and $\beta_i=0$ for $i\ne i_1,i_2$, $(\beta_i)_{i\in B}$ is represented as 
\[
(\beta_i)_{i\in B}=(p,q)_{i_1,i_2}.
\]
Using the above notations, we define moments related to $g_N$ as follows:
\[
M_0:=\int_{\mathbb{R}^{N-1}}g_N(\bmz',0)\,d\bmz', \quad  
M_{(\beta_i)_{i\in B}}^{\ell,m}:=\int_{\mathbb{R}^{N-1}}
\bigl\{(-\partial_{z_N}^2)^\ell(-\Delta_{\sbmz})^mg_N(\bmz',0)\bigr\}\,(\bmz')^{\beta}\,d\bmz'
\]
for $\ell,m\in\mathbb{Z}_+$ and $\beta=(\beta_1,\cdots,\beta_{N-1})\in \mathbb{Z}_+^{N-1}$.
We calculate the explicit values of some moments. 

\begin{lem}\label{lem:moment}
Let $\ell,m\in\mathbb{Z}_+$ and $\beta=(\beta_1,\cdots,\beta_{N-1})\in\mathbb{Z}_+^{N-1}$
with $|\beta|:=\beta_1+\cdots+\beta_{N-1}\le 6 $. Then, the following (i), (ii) and (iii) hold:
\begin{list}{}{\topsep=0.1cm\itemsep=0cm\leftmargin=0.3cm}
\item[\,(i)]
$M_0=c_1L_0$. 
\item[(ii)]
If $|\beta|$ is odd, $M_{(\beta_i)_{i\in B}}^{\ell,m}=0$. 
\item[(iii)]
For $i_0,i_1,i_2,i_3\in\Lambda$ with $i_1<i_2<i_3$,
\begin{align*}
&M_{(\beta_i)_{i\in B}}^{0,0}
=\left\{\begin{array}{cl} 
4c_1L_2&((\beta_i)_{i\in B}=(2)_{i_0}), \\[0.05cm]
-12c_1L_0&((\beta_i)_{i\in B}=(4)_{i_0}), \\[0.05cm]
-4c_1L_0&((\beta_i)_{i\in B}=(2,2)_{i_1,i_2}), 
\end{array}\right. \\
&M_{(\beta_i)_{i\in B}}^{1,0}
=\left\{\begin{array}{cl} 
-60c_1L_0&((\beta_i)_{i\in B}=(6)_{i_0}), \\
-12c_1L_0&((\beta_i)_{i\in B}=(4,2)_{i_1,i_2}), \\
-4c_1L_0&((\beta_i)_{i\in B}=(2,2,2)_{i_1,i_2,i_3}),
\end{array}\right.  \\
&M_{(2)_{i_0}}^{0,1}=-c_1L_0.
\end{align*}
\end{list}
\end{lem}
\noindent
{\bf Proof.} 
{\it Step 1.} We first prove (i). 
It follows from Lemma \ref{lem:g_xN=0} that for $\beta\in\mathbb{Z}_+^{N-1}$
\begin{align*}
&\int_{\mathbb{R}^{N-1}}
\bigl\{(-\partial_{z_N}^2)^\ell(-\Delta_{\sbmz})^mg_N(\bmz',0)\bigr\}\,(\bmz')^{\beta}\,d\bmz' \\
&=c_1\sum_{j=0}^m\binom{m}{j}\int_{\mathbb{R}^{N-1}}\mathcal{F}^{-1}_{N-1}
\biggl[\,\sum_{k=0}^\infty
\frac{(-2)^k}{k!}L_{2(k+j+\ell)}|\,\cdot\,|^{2(k+m-j)}e^{-|\,\cdot\,|^4}\,\biggr](\bmz')\,(\bmz')^{\beta}\,d\bmz'.
\end{align*}
Setting $(\ell,m)=(0,0)$ and $\beta=(0,\cdots,0)\in\mathbb{Z}_+^{N-1}$, we have
\begin{align*}
M_0
=&\,c_1\mathcal{F}_{N-1}\biggl[\mathcal{F}^{-1}_{N-1}\biggl[\,
\sum_{k=0}^\infty\frac{(-2)^{k}}{m!}L_{2k}|\,\cdot\,|^{2k}e^{-|\,\cdot\,|^4}\,\biggr]\biggr]
(\bmzero') \\
=&\,c_1\left.\sum_{k=0}^\infty\frac{(-2)^{k}}{k!}L_{2k}|\bmxi'|^{2k}
e^{-|\sbmxi'|^4}\right|_{\sbmxi'=\sbmzero'}
=c_1L_0.
\end{align*}

{\it Step 2.} We derive the precise form of $M_{(\beta_i)_{i\in B}}^{\ell,m}$ and prove (ii). It is observed that 
\begin{align*}
&\int_{\mathbb{R}^{N-1}}\mathcal{F}^{-1}_{N-1}
\biggl[\,\sum_{k=0}^\infty
\frac{(-2)^k}{k!}L_{2(k+j+\ell)}|\,\cdot\,|^{2(k+m-j)}e^{-|\,\cdot\,|^4}\,\biggr](\bmz')\,(\bmz')^{\beta}\,d\bmz' \\
&=(-\textbf{i})^{-|\beta|}\int_{\mathbb{R}^{N-1}}\mathcal{F}^{-1}_{N-1}
\biggl[\,\sum_{k=0}^\infty\frac{(-2)^k}{k!}L_{2(k+j+\ell)}
|\,\cdot\,|^{2(k+m-j)}e^{-|\,\cdot\,|^4}\,\biggr](\bmz')\,(-\textbf{i}\bmz')^{\beta}\,d\bmz' \\
&=\textbf{i}^{|\beta|}\int_{\mathbb{R}^{N-1}}\mathcal{F}^{-1}_{N-1}
\biggl[\,\sum_{k=0}^\infty\frac{(-2)^k}{k!}L_{2(k+j+\ell)}D_{\sbmxi'}^\beta
\bigl(|\bmxi'|^{2(k+m-j)}e^{-|\sbmxi'|^4}\bigr)\,\biggr](\bmz')\,\,d\bmz' \\
&=\textbf{i}^{|\beta|}\sum_{k=0}^\infty
\frac{(-2)^k}{k!}L_{2(k+j+\ell)}D_{\sbmxi'}^\beta\bigl(|\bmxi'|^{2(k+m-j)}e^{-|\sbmxi'|^4}\bigr)(\bmzero').
\end{align*}
This implies that 
\[
M_{(\beta_i)_{i\in B}}^{\ell,m}
=c_1\,\textbf{i}^{|\beta|}\sum_{j=0}^m\binom{m}{j}\sum_{k=0}^\infty\frac{(-2)^k}{k!}L_{2(k+j+\ell)}
D_{\sbmxi'}^\beta\bigl(|\bmxi'|^{2(k+m-j)}e^{-|\sbmxi'|^4}\bigr)(\bmzero').
\]
In particular, we see that 
\begin{align}
\label{moment00}
M_{(\beta_i)_{i\in B}}^{0,0}
=&\,c_1\,\textbf{i}^{|\beta|}\sum_{k=0}^\infty\frac{(-2)^k}{k!}L_{2k}D_{\sbmxi'}^\beta
\bigl(|\bmxi'|^{2k}e^{-|\sbmxi'|^4}\bigr)(\bmzero'), \\
\label{moment10}
M_{(\beta_i)_{i\in B}}^{1,0}
=&\,c_1\,\textbf{i}^{|\beta|}\sum_{k=0}^\infty\frac{(-2)^k}{k!}L_{2(k+1)}D_{\sbmxi'}^\beta
\bigl(|\bmxi'|^{2k}e^{-|\sbmxi'|^4}\bigr)(\bmzero'), \\
\label{moment01}
M_{(\beta_i)_{i\in B}}^{0,1}
=&\,c_1\,\textbf{i}^{|\beta|}\biggl\{
\sum_{k=0}^\infty\frac{(-2)^k}{k!}L_{2k}D_{\sbmxi'}^\beta\bigl(|\bmxi'|^{2(k+1)}e^{-|\sbmxi'|^4}\bigr)(\bmzero') \\
&\,\hspace*{35pt}
+\sum_{k=0}^\infty\frac{(-2)^k}{k!}L_{2(k+1)}D_{\sbmxi'}^\beta\bigl(|\bmxi'|^{2k}e^{-|\sbmxi'|^4}\bigr)(\bmzero')\biggr\}. 
\nonumber
\end{align}
Since $M_{(\beta_i)_{i\in B}}^{\ell,m}$ is real valued by its definition and $ \textbf{i}^{|\be|}=\pm \textbf{i} $ 
if $ |\be| $ is odd, we easily obtain (ii). 

{\it Step 3.} We show that (iii) holds. 
By virtue of the Leibniz rule, we have
\[
D_{\sbmxi'}^{\beta}\bigl(|\bmxi'|^{2(k+j)}e^{-|\sbmxi'|^4}\bigr)
=\sum_{\gamma\le\beta}\binom{\beta}{\gamma}\bigl(D^{\gamma}|\bmxi'|^{2(k+j)}\bigr)
\bigl(D^{\beta-\gamma}e^{-|\sbmxi'|^4}\bigr)
\]
for $j=0,1$, where $\gamma=(\gamma_1,\cdots,\gamma_{N-1})\in\mathbb{Z}_+^{N-1}$ and
\[
\binom{\beta}{\gamma}=\binom{\beta_1}{\gamma_1}\cdots\binom{\beta_{N-1}}{\gamma_{N-1}}, \quad
\binom{\beta_i}{\gamma_i}=\frac{\beta_i!}{\gamma_i!(\beta_i-\gamma_i)!}.
\]
Note that $\gamma\le\beta$ is defined as $\gamma_i\le\beta_i\,(i=1,\cdots,N-1)$. Since 
\[
|\bmxi'|^{2(k+j)}
=(\xi_1^2+\cdots+\xi_{N-1}^2)^{k+j}
=\sum_{d_1+\cdots+d_{N-1}=k+j}\binom{k+j}{d_1,\cdots,d_{N-1}}\xi_1^{2d_1}\cdots\xi_{N-1}^{2d_{N-1}},
\]
we see that 
\begin{align*}
&D_{\sbmxi'}^{\gamma}(|\bmxi'|^{2(k+j)})(\bmzero') \\
&=\sum_{d_1+\cdots+d_{N-1}=k+j}\binom{k+j}{d_1,\cdots,d_{N-1}}
\bigl(\partial_{\xi_1}^{\gamma_1}\xi_1^{2d_1}\bigr)
\cdots\bigl(\partial_{\xi_{N-1}}^{\gamma_{N-1}}\xi_{N-1}^{2d_{N-1}}\bigr)\biggr|_{\sbmxi'=\sbmzero'} \\
&=\left\{\begin{array}{cl}
\displaystyle\binom{k+j}{\gamma_1/2,\cdots,\gamma_{N-1}/2}\gamma_1!\cdots\gamma_{N-1}!
&(\gamma_i \text{ is even},\ \gamma_1+\cdots+\gamma_{N-1}=2(k+j)), \\[0.4cm]
0&\text{(otherwise)} .
\end{array}\right.
\end{align*}
This implies that 
\begin{align}\label{beta_derivative_2m}
&D_{\sbmxi'}^{\beta}\bigl(|\bmxi'|^{2(k+j)}e^{-|\sbmxi'|^4}\bigr)(\bmzero') \\
&=\sum_{\substack{2\sigma_1\le\beta_1,\cdots,2\sigma_{N-1}\le\beta_{N-1},\\ \sigma_1+\cdots+\sigma_{N-1}=k+j}}
\binom{\beta_1}{2\sigma_1}\cdots\binom{\beta_{N-1}}{2\sigma_{N-1}}
\binom{k+j}{\sigma_1,\cdots,\sigma_{N-1}} \nonumber \\
&\hspace*{110pt}(2\sigma_1)!\cdots(2\sigma_{N-1})!\bigl(\partial_{\xi_1}^{\beta_1-2\sigma_1}
\cdots\partial_{\xi_{N-1}}^{\beta_{N-1}-2\sigma_{N-1}}
e^{-|\sbmxi'|^4}\bigr)(\bmzero') \nonumber \\[0.1cm]
&=\sum_{\substack{2\sigma_1\le\beta_1,\cdots,2\sigma_{N-1}\le\beta_{N-1},\\ \sigma_1+\cdots
+\sigma_{N-1}=k+j}}\frac{\beta_1!\cdots\beta_{N-1}!}
{(\beta_1-2\sigma_1)!\cdots(\beta_{N-1}-2\sigma_{N-1})!}\binom{k+j}{\sigma_1,\cdots,\sigma_{N-1}} \nonumber \\
&\hspace*{110pt}\bigl(\partial_{\xi_1}^{\beta_1-2\sigma_1}\cdots\partial_{\xi_{N-1}}^{\beta_{N-1}
-2\sigma_{N-1}}e^{-|\sbmxi'|^4}\bigr)(\bmzero'). \nonumber
\end{align}
Here it follows that 
\begin{align*}
&\partial_{\xi_i}e^{-|\sbmxi'|^4}\bigr|_{\sbmxi'=\sbmzero}
=\partial_{\xi_{i_1}}\partial_{\xi_{i_2}}e^{-|\sbmxi'|^4}\bigr|_{\sbmxi'=\sbmzero}
=\partial_{\xi_{i_1}}\partial_{\xi_{i_2}}\partial_{\xi_{i_3}}e^{-|\sbmxi'|^4}\bigr|_{\sbmxi'=\sbmzero}=0, \\
&\partial_{\xi_{i_1}}\partial_{\xi_{i_2}}\partial_{\xi_{i_3}}\partial_{\xi_{i_4}}e^{-|\sbmxi'|^4}\bigr|_{\sbmxi'=\sbmzero}
=-8(\delta_{i_1i_2}\delta_{i_3i_4}+\delta_{i_1i_3}\delta_{i_2i_4}+\delta_{i_1i_4}\delta_{i_2i_3}), \\
&\partial_{\xi_{i_1}}\partial_{\xi_{i_2}}\partial_{\xi_{i_3}}\partial_{\xi_{i_4}}
\partial_{\xi_{i_5}}e^{-|\sbmxi'|^4}\bigr|_{\sbmxi'=\sbmzero}
=\partial_{\xi_{i_1}}\partial_{\xi_{i_2}}\partial_{\xi_{i_3}}\partial_{\xi_{i_4}}
\partial_{\xi_{i_5}}\partial_{\xi_{i_6}}e^{-|\sbmxi'|^2}\bigr|_{\sbmxi'=\sbmzero}=0.
\end{align*}
Therefore, we obtain 
\begin{align}\label{beta_derivative_4}
&\bigl(\partial_{\xi_1}^{\beta_1-2\sigma_1}\cdots\partial_{\xi_{N-1}}^{\beta_{N-1}-2\sigma_{N-1}}
e^{-|\sbmxi'|^4}\bigr)(\bmzero') \\[0.05cm] 
&=\left\{\begin{array}{cl} 
\,\ \ 1&(\beta_i-2\sigma_i=0\ (i=1,\cdots,N-1)), \\
-4!&(\beta_{i_0}-2\sigma_{i_0}=4,\,\beta_i-2\sigma_i=0\,(i\ne i_0)), \\
-8&(\beta_i-2\sigma_i=2\,(i= i_1,i_2),\,\beta_i-2\sigma_i=0\,(i\ne i_1,i_2)), \\
\,\ \ 0&(\mbox{otherwise})
\end{array}\right. \nonumber
\end{align}
for $|\beta|\le6$ and $i_0,i_1,i_2,i_3\in\Lambda$ with $i_1<i_2<i_3$. 
Thus \eqref{beta_derivative_2m} and \eqref{beta_derivative_4} imply that 
\begin{align}\label{derivative_2k_zero}
&D_{\sbmxi'}^{\beta}\bigl(|\bmxi'|^{2(k+j)}e^{-|\sbmxi'|^4}\bigr)(\bmzero') \\
&=\left\{\begin{array}{cl} 
\,\ \ 2&((\beta_i)_{i\in B}=(2)_{i_0},\,(k,j)=(1,0),(0,1)), \\
-4!&((\beta_i)_{i\in B}=(4)_{i_0},\,(k,j)=(0,0)), \\
\,\ \ 4!&((\beta_i)_{i\in B}=(4)_{i_0},\,(k,j)=(2,0),(1,1)), \\
-8&((\beta_i)_{i\in B}=(2,2)_{i_1,i_2},\,(k,j)=(0,0)), \\
\,\ \ 8&((\beta_i)_{i\in B}=(2,2)_{i_1,i_2},\,(k,j)=(2,0),(1,1)), \\
-6!&((\beta_i)_{i\in B}=(6)_{i_0},\,(k,j)=(1,0),(0,1)) \\
\,\ \ 6!&((\beta_i)_{i\in B}=(6)_{i_0},\,(k,j)=(3,0),(2,1)), \\
-6\cdot4!&((\beta_i)_{i\in B}=(4,2)_{i_1,i_2},\,(k,j)=(1,0),(0,1)), \\
\,\ \ 6\cdot4!&((\beta_i)_{i\in B}=(4,2)_{i_1,i_2},\,(k,j)=(3,0),(2,1)), \\
-2\cdot4!&((\beta_i)_{i\in B}=(2,2,2)_{i_1,i_2,i_3},\,(k,j)=(1,0),(0,1)), \\
\,\ \ 2\cdot4!&((\beta_i)_{i\in B}=(2,2,2)_{i_1,i_2,i_3},\,(k,j)=(3,0),(2,1)), \\
\,\ \ 0&(\mbox{otherwise})
\end{array}\right. \nonumber
\end{align}
for $|\beta|\le6$. 
Furthermore, by  \eqref{moment00}, \eqref{moment10}, \eqref{moment01} 
and \eqref{derivative_2k_zero}, we have
\begin{align*}
&M_{(\beta_i)_{i\in B}}^{0,0}
=\left\{\begin{array}{cl} 
4c_1L_2&((\beta_i)_{i\in B}=(2)_{i_0}), \\[0.05cm]
-4!c_1(L_0-2L_4)&((\beta_i)_{i\in B}=(4)_{i_0}), \\[0.05cm]
-8c_1(L_0-2L_4)&((\beta_i)_{i\in B}=(2,2)_{i_1,i_2}), 
\end{array}\right. \\
&M_{(\beta_i)_{i\in B}}^{1,0}
=\left\{\begin{array}{cl} 
-6!c_1\biggl(2L_4-\dfrac8{3!}L_8\biggr)&((\beta_i)_{i\in B}=(6)_{i_0}), \\
-6\cdot4!c_1\biggl(2L_4-\dfrac8{3!}L_8\biggr)&((\beta_i)_{i\in B}=(4,2)_{i_1,i_2}), \\
-2\cdot4!c_1\biggl(2L_4-\dfrac8{3!}L_8\biggr)&((\beta_i)_{i\in B}=(2,2,)_{i_1,i_2,i_3}),
\end{array}\right.  \\
&M_{(2)_{i_0}}^{0,1}
=-2c_1(L_0-2L_4).
\end{align*}
Since it follows from \eqref{cofficient-L} that
\[
L_4=\frac12\Gamma\biggl(\frac54\biggr)=\frac18\Gamma\biggl(\frac14\biggr)=\frac14L_0, \quad
L_8=\frac12\Gamma\biggl(\frac94\biggr)=\frac5{32}\Gamma\biggl(\frac14\biggr)=\frac5{16}L_0, 
\]
we are led to (iii). \qed

\bigskip
Combining Lemma \ref{lem:moment} with Theorem \ref{thm:asympt-exp_4w2}, 
we obtain the following precise asymptotic expansion of the solution to \eqref{4th_diffusion}. 

\begin{thm}\label{thm:asympt_expansion1}
Let $\gamma$ be a constant obtained in Lemma \ref{lem:series_finite} and take $v=Vt$ 
in Theorem \ref{thm:asympt-exp_4w2} for $V\in\mathbb{R}$. 
Then, an asymptotic expansion of a solution to \eqref{4th_diffusion} is given by
\begin{align}\label{asympt_expansion}
&\hspace*{-10pt}
u(\bmx+Vt\,\bmn(\bmx),t) \\
=&\,\frac12+c_1\Gamma\biggl(\frac34\biggr)H\,t^{1/4} \nonumber \\
&\,-\frac{c_1}2\Gamma\biggl(\frac14\biggr)\biggl\{V+\frac12\biggl(\Delta_gH+H|A|^2-\frac12H^3
+2\sum_{\substack{i_1,i_2,i_3\in\Lambda \\ i_1<i_2<i_3}}
\kappa_{i_1}\kappa_{i_2}\kappa_{i_3}-\lambda H\biggr)\biggr\}\,t^{3/4} \nonumber \\
&\,+O_{\p\om_0}(t) \nonumber
\end{align}
for any $\bmx\in\partial \Omega_0$, $\lambda\in\mathbb{R}$ and 
$t>0$ satisfying $|\lambda|t^{1/2}\le\gamma$ and $|V|t\le \delta_0$. 
\end{thm}

\noindent
{\bf Proof.}  
{\it Step 1.} Set 
\[
\Xi^{(0)}:=\int_{\mathbb{R}^{N-1}}\int_0^{\psi(\sbmz',Vt,t)}g_N(\bmz)\,dz_Nd\bmz', \quad
\Xi^{(1)}:=\int_{\mathbb{R}^{N-1}}\int_0^{\psi(\sbmz',Vt,t)}(-\Delta_{\sbmz})g_N(\bmz)\,dz_Nd\bmz'.
\]
Then we prove that
\begin{equation}
\label{eq0301}
    u(\bmx+Vt\,\bmn(\bmx),t)=\frac12+\Xi^{(0)}-\lambda t^{1/2}\,\Xi^{(1)}+O_{\p\om_0}(t).   
\end{equation}
According to \eqref{sol_expansion}, we see that
\begin{align*}
&\hspace*{-10pt}
u(\bmx+Vt\,\bmn(\bmx),t) \\
=&\,\frac12+\Xi^{(0)}
-\lambda t^{1/2}\Xi^{(1)} \\
&\,+\sum_{m=2}^\infty\frac{(-\lambda)^mt^{m/2}}{m!}
\int_{\mathbb{R}^{N-1}}\int_0^{\psi(\sbmz',Vt,t)}(-\Delta_{\sbmz})^mg_N(\bmz)\,dz_Nd\bmz'
+O_{\p\om_0}(e^{-\mu_\ast t^{-1/3}}). 
\end{align*}
%
By calculation similar to the proof of Lemma \ref{lem:series_finite}, 
we are able to find $C>0$ and $\gamma>0$ such that
\begin{align*}
&\left|\sum_{m=2}^\infty\frac{(-\lambda)^mt^{m/2}}{m!}
\int_{\mathbb{R}^{N-1}}\int_0^{\psi(\sbmz',v,t)}(-\Delta_{\sbmz})^mg_N(\bmz)\,dz_Nd\bmz'\right| \\
&\le\sum_{m=2}^\infty\frac{|\lambda|^mt^{m/2}}{m!}\int_{\mathbb{R}^N}|(-\Delta_{\sbmz})^mg_N(\bmz)|\,d\bmz
\le C\,|\lambda|^2\,t
\end{align*}
for $|\lambda|t^{1/2}\le\gamma$.  
This implies \eqref{eq0301}. 

{\it Step 2.} We derive the precise asymptotics for $ \Xi^{(0)} $ and $ \Xi^{(1)} $. 
Applying the Taylor's expansion of $g_N(\bmz)$ and $(-\Delta_{\sbmz})g_N(\bmz)$ 
with respect to the variable $z_N$ and taking account of 
\begin{equation}\label{Dg_zN=0}
\partial_{z_N}g_N(\bmz',0)=\partial_{z_N}^3g_N(\bmz',0)=0, \quad \partial_{z_N}(-\Delta_{\sbmz})g_N(\bmz',0)=0,
\end{equation}
which is obtained with the help of \eqref{1ode-phi_N}, we have
\begin{align*}
&g_N(\bmz)=g_N(\bmz',0)+\frac12\partial_{z_N}^2g_N(\bmz',0)z_N^2
+\frac1{4!}\partial_{z_N}^4g_N(\bmz',\theta_0z_N)z_N^4, \\
&(-\Delta_{\sbmz})g_N(\bmz)=(-\Delta_{\sbmz})g_N(\bmz',0)
+\frac12\partial_{z_N}^2(-\Delta_{\sbmz})g_N(\bmz',\theta_1z_N)z_N^2
\end{align*}
for some $\theta_0,\theta_1\in(0,1)$. These imply that 
\begin{align*}
\Xi^{(0)}=&\,\int_{\mathbb{R}^{N-1}}g_N(\bmz',0)\,\psi(\bmz',Vt,t)\,d\bmz'
-\frac16\int_{\mathbb{R}^{N-1}}(-\partial_{z_N}^2)g_N(\bmz',0)\bigl\{\psi(\bmz',Vt,t)\bigr\}^3\,d\bmz' \\
&\,+\frac1{4!}\int_{\mathbb{R}^{N-1}}\int_0^{\psi(\sbmz',Vt,t)}
(-\partial_{z_N}^2)^2g_N(\bmz',\theta_0z_N)z_N^4\,dz_Nd\bmz', \\
\Xi^{(1)}=&\,\int_{\mathbb{R}^{N-1}}(-\Delta_{\sbmz})g_N(\bmz',0)\,\psi(\bmz',Vt,t)d\bmz' \\
&\,-\frac12\int_{\mathbb{R}^{N-1}}\int_0^{\psi(\sbmz',Vt,t)}
(-\partial_{z_N}^2)(-\Delta_{\sbmz})g_N(\bmz',\theta_1z_N)z_N^2\,dz_Nd\bmz'. 
\end{align*}
Using Corollary \ref{cor:Cui-thm}, we observe that 
\begin{align*}
&\left|\int_{\mathbb{R}^{N-1}}\int_0^{\psi(\sbmz',Vt,t)}
(-\partial_{z_N}^2)^2g_N(\bmz',\theta_0z_N)z_N^4\,dz_Nd\bmz'\right|=O_{\p\om_0}(t^{5/4}), \\
&\left|\int_{\mathbb{R}^{N-1}}\int_0^{\psi(\sbmz',Vt,t)}
(-\partial_{z_N}^2)(-\Delta_{\sbmz})g_N(\bmz',\theta_1z_N)z_N^2\,dz_Nd\bmz'\right|
=O_{\p\om_0}(t^{3/4}).
\end{align*}
Then it follows from \eqref{surface_expansion}, \eqref{kappa1}, \eqref{kappa3} and Lemma \ref{lem:moment} that
\begin{align*}
\Xi^{(0)}=&\,\frac12\int_{\mathbb{R}^{N-1}}g_N(\bmz',0)\sum_{i_0\in\Lambda}z_{i_0}^2\,\kappa_{i_0}\,d\bmz'\,\,t^{1/4}
-\int_{\mathbb{R}^{N-1}}g_N(\bmz',0)V\,d\bmz'\,\,t^{3/4} \\
&\,+\frac1{24}\int_{\mathbb{R}^{N-1}}g_N(\bmz',0)\biggl(\sum_{i_0\in\Lambda}z_{i_0}^4\partial_{z_{i_0}}^4f
+6\sum_{\substack{i_1,i_2\in\Lambda \\ i_1<i_2}}
z_{i_1}^2z_{i_2}^2\partial_{z_{i_1}}^2\partial_{z_{i_2}}^2f\biggr)\,d\bmz'\,\,t^{3/4} \\
&\,-\frac1{48}\int_{\mathbb{R}^{N-1}}(-\partial_{z_N}^2)g_N(\bmz',0)
\biggl(\sum_{i_0\in\Lambda}z_{i_0}^6\kappa_{i_0}^3
+3\sum_{\substack{i_1,i_2\in\Lambda \\ i_1\ne i_2}}z_{i_1}^4z_{i_2}^2\kappa_{i_1}^2\kappa_{i_2} \\
&\hspace*{150pt}
+6\sum_{\substack{i_1,i_2,i_3\in\Lambda \\ i_1<i_2<i_3}}
z_{i_1}^2z_{i_2}^2z_{i_3}^2\kappa_{i_1}\kappa_{i_2}\kappa_{i_3}\biggr)\,d\bmz'\,\,t^{3/4}
+O_{\p\om_0}(t) \\
=&\,\frac12\sum_{i_0\in\Lambda}M_{(2)_{i_0}}^{0,0}\,\kappa_{i_0}\,t^{1/4}-M_0\,V\,t^{3/4} \\
&\,+\frac1{24}\biggl(\sum_{i_0\in\Lambda}M_{(4)_{i_0}}^{0,0}\partial_{z_{i_0}}^4f
+6\sum_{\substack{i_1,i_2\in\Lambda \\ i_1<i_2}}M_{(2,2)_{i_1,i_2}}^{0,0}
\partial_{z_{i_1}}^2\partial_{z_{i_2}}^2f\biggr)\,t^{3/4} \\
&\,-\frac1{48}\biggl(\sum_{i_0\in\Lambda}M_{(6)_{i_0}}^{1,0}\kappa_{i_0}^3
+3\sum_{\substack{i_1,i_2\in\Lambda \\ i_1\ne i_2}}M_{(4,2)_{i_1,i_2}}^{1,0}\kappa_{i_1}^2\kappa_{i_2} \\
&\,\hspace*{40pt}
+6\sum_{\substack{i_1,i_2,i_3\in\Lambda \\ i_1<i_2<i_3}}
M_{(2,2,2)_{i_1,i_2,i_3}}^{1,0}\kappa_{i_1}\kappa_{i_2}\kappa_{i_3}
\biggr)\,t^{3/4}+O_{\p\om_0}(t) \\
=&\,2c_1L_2\sum_{i_0\in\Lambda}\kappa_{i_0}\,t^{1/4} \\
&\,-c_1L_0\biggl\{V+\frac12\biggl(\sum_{i_0\in\Lambda}\partial_{z_{i_0}}^4f
+\sum_{\substack{i_1,i_2\in\Lambda \\ i_1<i_2}}\partial_{z_{i_1}}^2\partial_{z_{i_2}}^2f
-\frac52\sum_{i_0\in\Lambda}\kappa_{i_0}^3-\frac32\sum_{\substack{i_1,i_2\in\Lambda \\ i_1\ne i_2}}\kappa_{i_1}^2\kappa_{i_2} \\
&\,\hspace*{85pt}
-\sum_{\substack{i_1,i_2,i_3\in\Lambda \\ i_1<i_2<i_3}}\kappa_{i_1}\kappa_{i_2}\kappa_{i_3}\biggr)\biggr\}\,t^{3/4}+O_{\p\om_0}(t), \\
\Xi^{(1)}=&\,\frac12\int_{\mathbb{R}^{N-1}}
(-\Delta_{\sbmz})g_N(\bmz',0)\sum_{i_0\in\Lambda}z_{i_0}^2\,\kappa_{i_0}\,d\bmz'\,\,t^{1/4}+O_{\p\om_0}(t^{3/4}) \\
=&\,-\frac12c_1L_0\sum_{i_0\in\Lambda}\kappa_{i_0}\,t^{1/4}+O_{\p\om_0}(t^{3/4}).
\end{align*}

{\it Step 3.} From the results in Step 1 and Step 2, we prove (\ref{asympt_expansion}). 
Taking account of $L_0=\Gamma(1/4)/2$ and $L_2=\Gamma(3/4)/2$, we obtain
\begin{align*}
&\hspace*{-10pt}
u(\bmx+Vt\,\bmn(\bmx),t) \\
=&\,\frac12+\Xi^{(0)}-\lambda t^{1/2}\,\Xi^{(1)}+O_{\p\om_0}(t) \\
=&\,\frac12+c_1\Gamma\biggl(\frac34\biggr)\sum_{i_0\in\Lambda}\kappa_{i_0}\,t^{1/4} \\
&\,-\frac{c_1}2\Gamma\biggl(\frac14\biggr)\biggl\{V+\frac12\biggl(\sum_{i_0\in\Lambda}\partial_{z_{i_0}}^4f
+2\sum_{\substack{i_1,i_2\in\Lambda \\ i_1<i_2}}\partial_{z_{i_1}}^2\partial_{z_{i_2}}^2f
-\frac52\sum_{i_0\in\Lambda}\kappa_{i_0}^3-\frac32\sum_{\substack{i_1,i_2\in\Lambda \\ i_1\ne i_2}}\kappa_{i_1}^2\kappa_{i_2} \\
&\hspace*{110pt}
-\sum_{\substack{i_1,i_2,i_3\in\Lambda \\ i_1<i_2<i_3}}\kappa_{i_1}\kappa_{i_2}\kappa_{i_3}
-\lambda\sum_{i_0\in\Lambda}\kappa_{i_0}\biggr)\biggr\}\,t^{3/4}
+O_{\p\om_0}(t).
\end{align*}
Referring Appendix \ref{sec:mc} below, we see that 
\begin{align*}
&H=\sum_{i\in\Lambda}\kappa_i, \quad |A|^2=\sum_{i\in\Lambda}\kappa_i^2, \\
&\Delta_gH+H|A|^2-\frac12H^3 \\
&=\sum_{i_0\in\Lambda}\partial_{z_{i_0}}^4f
+2\sum_{\substack{i_1,i_2\in\Lambda \\ i_1<i_2}}\partial_{z_{i_1}}^2\partial_{z_{i_2}}^2f
-\frac52\sum_{i\in\Lambda}\kappa_i^3
-\frac32\sum_{\substack{i_1,i_2\in\Lambda \\ i_1\ne i_2}}\kappa_{i_1}^2\kappa_{i_2}
-3\sum_{\substack{i_1,i_2,i_3\in\Lambda \\ i_1<i_2<i_3}}
\kappa_{i_1}\kappa_{i_2}\kappa_{i_3}.   
\end{align*}
Consequently, we have \eqref{asympt_expansion}.  \qed

\section{A thresholding algorithm to the flow by (\ref{geometric-flow})}
In this section,  we first introduce a threshold function based on the asymptotic expansion obtained in Theorem \ref{thm:asympt_expansion1} 
and propose the thresholding algorithm. After that, we justify its algorithm. 
Throughout this section, we assume that $\Omega_0\subset\bR^N$ is a compact set and $\p\om_0$ is of class $C^5$.

\subsection{A thresholding algorithm}\label{subsec:proposed_algorithm}
Let $u$ be a solution to \eqref{4th_diffusion}. Set
\[
u_a(\bmx,t):=u(\bmx,a^4t)
\]
for $\bmx\in\mathbb{R}^N$ and $a>0$. Then, 
it follows from Theorem \ref{thm:asympt_expansion1} that for $V\in\mathbb{R}$ 
\begin{align*}
&\hspace*{-10pt}
u_a(\bmx+Vt\,\bmn(\bmx),t) \\
=&\,\frac12+c_1\Gamma\biggl(\frac34\biggr)aH\,t^{1/4} \\
&\,-\frac{c_1}2\Gamma\biggl(\frac14\biggr)
\biggl\{\frac{V}{a}+\frac{a^3}2\biggl(\Delta_gH+H|A|^2-\frac12H^3
+2\sum_{\substack{i_1,i_2,i_3\in\Lambda \\ i_1<i_2<i_3}}\kappa_{i_1}\kappa_{i_2}\kappa_{i_3}-\lambda H\biggr)\biggr\}\,t^{3/4} \\
&\,+O_{\p\om_0}(t).  
\end{align*}
Hence we see that 
%
%
\begin{align*}
&\hspace*{-10pt}
u_{3a}(\bmx+Vt\,\bmn(\bmx),t)-3u_{2a}(\bmx+Vt\,\bmn(\bmx),t)+3u_a(\bmx+Vt\,\bmn(\bmx),t)-\frac12 \\
=&\,-\frac{11c_1}{12a}\Gamma\biggl(\frac14\biggr) \\
&\,\cdot\biggl\{V+\frac{18a^4}{11}\biggl(\Delta_gH+H|A|^2-\frac12H^3
+2\sum_{\substack{i_1,i_2,i_3\in\Lambda \\ i_1<i_2<i_3}}\kappa_{i_1}\kappa_{i_2}\kappa_{i_3}-\lambda H\biggr)\biggr\}\,t^{3/4} \\
&\,+O_{\p\om_0}(t).
\end{align*}
Choosing $a>0$ such that $18a^4/11=1$, we are led to 
\begin{align}
\label{esti040101}
&u_{3a}(\bmx+Vt\,\bmn(\bmx),t)-3u_{2a}(\bmx+Vt\,\bmn(\bmx),t)+3u_a(\bmx+Vt\,\bmn(\bmx),t)-\frac12 \\
\nonumber
&=-\frac{11c_1}{12a}\Gamma\biggl(\frac14\biggr)\biggl(V+\Delta_gH+H|A|^2-\frac12H^3
+2\sum_{\substack{i_1,i_2,i_3\in\Lambda \\ i_1<i_2<i_3}}\kappa_{i_1}\kappa_{i_2}\kappa_{i_3}-\lambda H\biggr)\,t^{3/4} \\
\nonumber
& \quad +O_{\p\om_0}(t).
\end{align}
From the above observation, let us introduce a threshold function and a new set generated by its function. 
Define a threshold function $U(\bmx,t)$ as 
\begin{equation}\label{threshold-func}
U(\bmx,t):=u_{3a}(\bmx,t)-3u_{2a}(\bmx,t)+3u_a(\bmx,t)\quad\text{}
\end{equation}
for $a>0$ satifying $18a^4/11=1$ and set 
\begin{equation}\label{set_Omega_t}
\Omega(t):=\Bigl\{\bmx\in\mathbb{R}^N\,\Big|\,U(\bmx,t)\ge\frac12\Bigr\}.  
\end{equation}
For any $ \bmx\in\p\om_0 $ and small $ t>0 $, we define $ V=V(\bmx,t) $ by 
\begin{equation}
\label{eq:def-V}
   \bmx+V(\bmx,t)t\,\bmn(\bmx)\in\p\Omega(t).  
\end{equation}
Then, setting $\bmy(\bmx,t):=\,\bmx+V(\bmx,t)t\,\bmn(\bmx) $, we notice that 
\[
    V(\bmx,t)=-\frac{d(\bmy(\bmx,t),\p\om_0)}{t}, \quad 
    |d(\bmy(\bmx,t),\p\om_0)|=|\bmy(\bmx,t)-\bmx|
\] 
for all $ \bmx\in\p\om_0 $ and small $ t>0 $.  Hence we are able to regard that $ V $ is the outward normal velocity 
from $ \p\om_0 $ to $ \p\Omega(t) $.  Here $ d(\bmy,\p\om_0) $ is defined by (\ref{signed-distance}).  

We assume that 
\begin{equation}
\label{eq:V-bounded}
	\|V\|_{L^\infty(\p\om_0\times(0, t_0))}<\infty\quad\text{for some}\,\ t_0>0,
\end{equation}
and $U(\bmx+V(\bmx,t)t\,\bmn(\bmx),t)=1/2$.  Then, by \eqref{esti040101}
\[
V+\Delta_gH+H|A|^2-\frac12H^3
+2\sum_{\substack{i_1,i_2,i_3\in\Lambda \\ i_1<i_2<i_3}}\kappa_{i_1}\kappa_{i_2}\kappa_{i_3}-\lambda H=O_{\p\om_0}(t^{1/4}), 
\]
where $ g=g(\bmx) $ and $ \kappa_i=\kappa_i(\bmx) $ for $ \bmx\in\partial\Omega_0 $.  
This implies that $V$ must be closed to $-\nabla_{L^2}\mathcal{E}_{\lambda}^N(\partial\Omega_0)$. 
We emphasize that the assumption \eqref{eq:V-bounded} is actually valid as we show in the next subsection.

Based on the above argument, let us derive a thresholding algorithm to \eqref{geometric-flow}.
First we solve the initial value problem \eqref{4th_diffusion} for the initial function $\chi_{\Omega_0}(\bmx)$ 
and let $u^0$ be the corresponding solution. 
Define a threshold function $U^0(\bmx,t)$ as \eqref{threshold-func} and a set $\Omega^0(t)$ as \eqref{set_Omega_t}. 
Fix a time step $h>0$ and define $\Omega_1:=\om^0(h) $. 
As the second step, we solve the problem \eqref{4th_diffusion} with $\Omega_1$ replacing $\Omega_0$, 
define $\Omega^1(t) $ as the set of \eqref{set_Omega_t} with $u$ replaced by the solution to the problem \eqref{4th_diffusion} 
with the new initial function $\chi_{\Omega_1}(\bmx)$.  
Repeat this procedure to obtain a sequence $\{\Omega_k\}_{k\in\mathbb{Z}_+}$ of compact sets in $\mathbb{R}^N$. 
Then, setting 
\[
\widehat{\Omega}^h(t)=\Omega_k\,\,\ \text{for}\,\ t\in[kh,(k+1)h)\,\ (k\in\mathbb{Z}_+)
\]
and letting $h\to0$, we can expect that at least formally there is a limit flow $\{\partial\widehat\Omega(t)\}_{t\ge0}$ of 
$ \{\p\widehat{\om}^h(t)\}_{t\geq 0} $ as $ h\to 0 $, whose boundary moves by \eqref{geometric-flow} with 
$ \widehat{\G}(t)=\p\widehat{\om}(t) $. 
Indeed, since $\widehat{V}^h(\bmx,t )$, given by 
\begin{align*}
 & \bmx+\widehat{V}^h(\bmx,t )h\bmn(\bmx)\in\p\om_{k+1}\,\,\ \text{for}\,\ t\in[kh,(k+1)h)\ (k\in\mathbb{Z}_+)\ 
   \text{and}\,\ \bmx\in\p\om_k,
\end{align*}
can be regarded as a normal velocity of $\p\widehat{\om}^h(t)$, the above observation implies that the limit flow of 
$\{\p\widehat{\om}^h(t)\}_{t\geq 0} $ formally moves by \eqref{geometric-flow}. 

\begin{rem}\label{rem:difference}
Set
\begin{align*}
&E_\alpha(\bmx,t):=\frac1{(t^\alpha)^N}\rho\left(\frac{|\bmx|}{t^\alpha}\right), \quad
\rho(s):=\frac1{(4\pi)^{N/2}}\exp\left(-\frac{s^2}4\right), \\
&\psi_\alpha(\bmz',v,t):=t^{-\alpha}\{-v+f(t^\alpha\bmz')\}.
\end{align*}
Then the expansion of $\psi_\alpha(\bmz',Vt,t)$ is given by 
\begin{align*}
\psi_\alpha(\bmz',Vt,t)
=&\,-Vt^{1-\alpha}+\frac{\langle\bmz',\nabla_{\sbmz'}\rangle_{N-1}^2f(\bmzero')}2t^\alpha
+\frac{\langle\bmz',\nabla_{\sbmz'}\rangle_{N-1}^3f(\bmzero')}6t^{2\alpha} \\
&\,+\frac{\langle\bmz',\nabla_{\sbmz'}\rangle_{N-1}^4f(\bmzero')}{24}t^{3\alpha}
+\frac{\langle\bmz',\nabla_{\sbmz'}\rangle_{N-1}^5f(\theta t^\alpha\bmz')}{120}t^{4\alpha}
\end{align*}
for some $\theta\in(0,1)$. 
Let $\nabla_{L^2}\mathcal{E}_0^N(\Gamma)$ be the $L^2$-gradient of $\mathcal{E}_0^N(\Gamma)$ 
given by \eqref{L2-gradient} with $\lambda=0$, that is, 
\[
\nabla_{L^2}\mathcal{E}_0^N(\Gamma)
=\Delta_gH+H|A|^2-\frac12H^3+2\sum_{\substack{i_1,i_2,i_3\in\Lambda \\ i_1<i_2<i_3}}
\kappa_{i_1}\kappa_{i_2}\kappa_{i_3}, 
\] 
and set  
\[
w_\alpha(\bmx,t):=(E_\alpha(\cdot,t)\ast\chi_{\Omega_0})(\bmx). 
\]
Since $E_{1/2}(\bmx,t)$ is the Gauss kernel, $w_{1/2}$ is a solution to 
\[
\left\{\begin{array}{l}
w_t=\Delta w\,\,\ \text{in}\ \mathbb{R}^N\times(0,\infty), \\
w(\bmx,0)=\chi_{\Omega_0}(\bmx)\,\,\ \text{in}\ \mathbb{R}^N.
\end{array}\right.
\]
Here the precise asymptotic expansion of $w_{1/2}(\bmx+Vt\,\bmn(\bmx),t)$ until the term of $t^{3/2}$ is represented as 
\begin{align*}
&\hspace*{-10pt}
w_{1/2}(\bmx+Vt\,\bmn(\bmx),t) \\
=&\,\frac12-c_1\sqrt{\pi}(V-H)\,t^{1/2} \nonumber \\
&\,+c_1\sqrt{\pi}\biggl\{\frac12\nabla_{L^2}\mathcal{E}_0^N(\Gamma)
+\frac14V(2|A|^2+H^2)-\frac14V^2H+\frac1{12}V^3\biggr\}\,t^{3/2}
+O_{\p\om_0}(t^2).
\end{align*}
We remark that the precise asymptotic expansion until the term of $t^{1/2}$ is obtained in Evans\,\cite{eva;93} and 
its expansion is the basis of the BMO algorithm for the mean curvature flow (cf. Bence--Merriman--Osher\,\cite{ben;mer;osh;92}). 
On the other hand, considering the asymptotic expansion of $w_{1/4}(\bmx+Vt\,\bmn(\bmx),t)$, 
we have
\begin{align}\label{w1/4-expansion}
&\hspace*{-10pt}
w_{1/4}(\bmx+Vt\,\bmn(\bmx),t) \\
=&\,\frac12+c_1\sqrt{\pi}H\,t^{1/4}-c_1\sqrt{\pi}\biggl(V-\frac12\nabla_{L^2}\mathcal{E}_0^N(\Gamma)\biggr)\,t^{3/4}
+O_{\p\om_0}(t). \nonumber
\end{align}
The difference between \eqref{w1/4-expansion} and our expansion 
\begin{align*}
&\hspace*{-10pt}
u(\bmx+Vt\,\bmn(\bmx),t) \\
=&\,\frac12+c_1\Gamma\biggl(\frac34\biggr)H\,t^{1/4}
-\frac{c_1}2\Gamma\biggl(\frac14\biggr)\biggl(V+\frac12\nabla_{L^2}\mathcal{E}_0^N(\Gamma)\biggr)\,t^{3/4}
+O_{\p\om_0}(t),
\end{align*}
which is given by \eqref{asympt_expansion} with $\lambda=0$, is the sign before $\nabla_{L^2}\mathcal{E}_0^N(\Gamma)$. 
This difference is due to the fact that $w_{1/4}$ satisfies the second order parabolic equation 
$w_t=(2\sqrt{t})^{-1}\Delta w$, whereas $u$ fulfills the fourth order parabolic equation $u_t=-\Delta^2u$. 
The thresholding algorithms of the Willmore flow based on \eqref{w1/4-expansion} are derived 
for $N=2$ in Esedo\=glu--Ruuth--Tsai\,\cite{ese;ruu;tsa;08} and for $N=3$ in 
Grzhibovskis--Heintz\,\cite{grz;hei;08}. 
Since the difference explained above is essentially related to the parabolicity, the same threshold function 
as ours, that is, $U$ given by \eqref{threshold-func} can not be chosen if $w_{1/4}$ is used. 
\end{rem}
\subsection{Properties of evolving sets} \label{sec:evolving-sets}
We give a justification of the argument in the previous subsection, more precisely, we prove that the assumption \eqref{eq:V-bounded} 
is actually valid. 
In order to do so, we prepare several propositions and lemmas. 
In the following, we use the notation defined as in \eqref{threshold-func}, \eqref{set_Omega_t} and \eqref{eq:def-V}. 

\begin{prop}\label{prop:inclusion}
There exists $K_\ast>0$ and $t_\ast>0$ such that 
\[
\p\Omega(t)\subset\{\bmx\in\mathbb{R}^N\,|\,|d(\bmx,\partial \Omega_0)|\le K_\ast t^{1/4}\}
\]
for $t\in(0,t_\ast)$, where $d(\,\cdot\,,\partial \Omega_0)$ is the signed distance function to 
$\partial \Omega_0$ given by \eqref{signed-distance}. 
\end{prop}
\noindent
{\bf Proof.} 
{\it Step 1}. 
Set $D_+(t,r):=\{\bmx\in\mathbb{R}^N\,|\,d(\bmx,\partial \Omega_0)>rt^{1/4}\}$ for $r>0$.  
We show that  there exists $ t_+>0 $ and $ K_*>0 $ such that $ D_+(t,K_*)\subset\Omega(t)$ 
for all $ t\in (0,t_+) $. Recalling the definition \eqref{signed-distance} of the signed distance function, 
we can find  $t_+=t_+(r)>0 $ satisfying  
$\emptyset\ne D_+(t,r)\subset \Omega_0$ for $t\in(0,t_+)$. Fix any $\bmx\in D_+(t,r)$. Taking account of  
$B_N(\bmx,rt^{1/4}/2)\subset \Omega_0$ 
for $t\in(0,t_+)$ and recalling \eqref{integralG=1}, we see that for any $a>0$
\begin{align*}
&|u_a(\bmx,t)-1| \\
&=|(G_{N,\lambda}(\,\cdot\,,a^4t)\ast\chi_{\Omega_0})(\bmx)-1| \\
&\le\left|\int_{B_N(\sbmx,rt^{1/4}/2)}G_{N,\lambda}(\bmx-\bmy,a^4t)\,d\bmy-1\right|
+\int_{\Omega_0\setminus B_N(\sbmx,rt^{1/4}/2)}|G_{N,\lambda}(\bmx-\bmy,a^4t)|\,d\bmy \\
&\le2\int_{\mathbb{R}^N\setminus B_N(\sbmx,rt^{1/4}/2)}|G_{N,\lambda}(\bmx-\bmy,a^4t)|\,d\bmy.
\end{align*}
Here it follows from Theorem \ref{thm:Cui} with $|\alpha|=0$ and $m=0$ that for $t\in(0,t_+)$ 
\begin{align*}
&\int_{\mathbb{R}^N\setminus B(\sbmx,rt^{1/4}/2)}|G_{N,\lambda}(\bmx-\bmy,a^4t)|\,d\bmy \\
&\le\frac{C_1e^{K|\lambda|^2t}}{a^Nt^{N/4}}
\int_{\mathbb{R}^N\setminus B(\sbmx,rt^{1/4}/2)}e^{-\mu(|\sbmx-\sbmy|^4/(a^4t))^{1/3}}\,d\bmy \\
&\le C_2\int_{\mathbb{R}^N\setminus B(\sbmzero,r/2a)}e^{-\mu|\sbmz|^{4/3}}\,d\bmz.
\end{align*}
Using the polar coordinate, we have
\[
\int_{\mathbb{R}^N\setminus B(\sbmzero,r/2a)}e^{-\mu|\sbmz|^{4/3}}\,d\bmz
=\omega_{N-1}\int_{r/2a}^\infty e^{-\mu\rho^{4/3}}\rho^{N-1}\,d\rho
\le C_3e^{-\mu_1(r/a)^{4/3}},
\]
where $C_3:=(3/4)(2/\mu)^{3N/4}\Gamma(3N/4)$ and $\mu_1:=\mu/2^{7/3}$. Thus it is seen that 
\[
|u_a(\bmx,t)-1|\le C_4\,e^{-\mu_1(r/a)^{4/3}}.
\]
By means of this inequality, we obtain
\[
U(\bmx,t)=u_{3a}(\bmx,t)-3u_{2a}(\bmx,t)+3u_a(\bmx,t)\ge1-7C_4\,e^{-\mu_1(r/3a)^{4/3}}
\]
for $t\in(0,t_+)$. Taking $r=K_\ast>0$ such that $7C_4\,e^{-\mu_1(K_\ast/3a)^{4/3}}<1/2$, we conclude that
\[
U(\bmx,t)>\frac12
\]
for $t\in(0,t_+(K_\ast))$. This implies that $\bmx\in{\rm int}\,\,\Omega(t)$, so that 
$D_+(t,K_\ast)\subset{\rm int}\,\,\Omega(t)$ for $t\in(0,t_+(K_\ast))$. 

{\it Step 2}. 
Set $D_-(t,r):=\{\bmx\in\mathbb{R}^N\,|\,d(\bmx,\partial \Omega_0)<-rt^{1/4}\}$ for $r>0$. 
We prove that there exists $t_->0$ such that $D_-(t,K_\ast)\subset\mathbb{R}^N\setminus\Omega(t)$ 
for all $t\in(0,t_-) $, where $K_\ast$ is a constant as in Step 1. Recalling \eqref{signed-distance} again, we are able to 
find a $ t_-=t_-(r) $ such that $\emptyset\ne D_-(t,r)\subset\mathbb{R}^N\setminus \Omega_0$ for $t\in(0,t_-)$. 
Fix any $\bmx\in D_-(t,r)$. Since $\Omega_0\subset\mathbb{R}^N\setminus B_N(\bmx,rt^{1/4}/2)$ for $t\in(0,t_-)$, 
it follows that for any $a>0$ and $t\in(0,t_-)$
\begin{align*}
|u_a(\bmx,t)|
=&\,|(G_{N,\lambda}(\,\cdot\,,a^4t)\ast\chi_{\Omega_0})(\bmx)| \\
\le&\,\int_{\mathbb{R}^N\setminus B(\sbmx,rt^{1/4}/2)}|G_{N,\lambda}(\bmx-\bmy,a^4t)|\,d\bmy \\
\le&\,C_4e^{-\mu_1(r/a)^{4/3}}.
\end{align*}
Using this inequality, we have
\[
U(\bmx,t)=u_{3a}(\bmx,t)-3u_{2a}(\bmx,t)+3u_a(\bmx,t)\le7C_4\,e^{-\mu_1(r/3a)^{4/3}}
\]
for $t\in(0,t_-)$. Taking $r=K_\ast>0$ as in Step 1, we obtain
\[
U(\bmx,t)<\frac12
\]
for $t\in(0,t_-(K_\ast))$. This implies that $\bmx\in\mathbb{R}^N\setminus\Omega(t)$, so that 
$D_-(t,K_\ast)\subset\mathbb{R}^N\setminus\Omega(t)$ for $t\in(0,t_-(K_\ast))$. 

{\it Step 3}. Define $t_\ast:=\min\{t_+(K_\ast),t_-(K_\ast)\}$. Then we see that 
\[
D_+(t,K_\ast)\subset{\rm int}\,\,\Omega(t), \quad 
D_-(t,K_\ast)\subset\mathbb{R}^N\setminus\Omega(t)
\]
for $t\in(0,t_\ast)$. This leads to the desired result. \qed

\bigskip\noindent
Furthermore, we are able to derive the following refinement of Proposition \ref{prop:inclusion}. 

\begin{prop}\label{prop:refinement}
For any $\varepsilon\in(0,1)$, there exists $t_{\ast,\varepsilon}>0$ such that 
\[
\p\Omega(t)\subset\{\bmx\in\mathbb{R}^N\,|\,|d(\bmx,\partial\Omega_0)|\le\varepsilon t^{1/4}\}
\]
for $t\in(0,t_{\ast,\varepsilon})$. 
\end{prop}

\noindent
In order to prove this proposition, we need several preparations. Let us consider the case where 
$ \lambda=0 $ and $ \Omega_0=\{\bmx\in\bR^N\;|\;x_N\le 0\}$.  
Taking account of $\bmn(\bmx)=\bme_N$ and $f\equiv0$ in this case, we see that 
for $\bmx\in\partial\Omega_0$, $v\in\mathbb{R}$ and $t>0$
\begin{align*}
u(\bmx+v\,\bmn(\bmx),t)=\frac12+\int_{\bR^{N-1}}\int_0^{-vt^{-1/4}}g_N(\bmz)\,d\bmz
=u(v\,\bme_N,t).  \nonumber
\end{align*}
Since $u(v\,\bme_N,t)$ is the solution to \eqref{4th_diffusion} with 
$ \lambda=0 $ and $ \Omega_0=\{\bmx\in\bR^N\;|\;x_N\le 0\}$, 
$u(v\,\bme_N,t)$ can be the solution to 
\begin{equation}\label{1D_4th_diffusion}
\left\{\begin{array}{l}
\widetilde{u}_t=-\partial_{y_N}^4\widetilde{u}
\,\,\ \text{in}\ \mathbb{R}\times(0,\infty), \\
\widetilde{u}(y_N,0)=\chi_{(-\infty,0]}(y_N)=\left\{\begin{array}{l} 1\,\,\ \text{in}\ (-\infty,0], \\ 
0\,\,\ \text{in}\ \mathbb{R}^N\setminus(-\infty,0]. \end{array}\right.
\end{array}\right.
\end{equation}
On the other hand, the solution $\widetilde{u}(y_N,t)$ to \eqref{1D_4th_diffusion} is represented as 
\[
\widetilde{u}(y_N,t)=\frac12+\int_0^{-y_Nt^{-1/4}}g_1(z)dz.  
\]
By the uniqueness of the solution to \eqref{1D_4th_diffusion}, we obtain
\[
u(v\,\bme_N,t)=\widetilde{u}(v,t)=\frac12+\int_0^{-vt^{-1/4}}g_1(z)dz. 
\]
As a result, it follows that for $\bmx\in\partial\Omega_0=\{\bmx\in\bR^N\;|\;x_N=0\}$
\begin{align}\label{sol_wholesp}
u(\bmx+v\,\bmn(\bmx),t) 
=\frac12+\int_{\bR^{N-1}}\int_0^{-vt^{-1/4}}g_N(\bmz)\,d\bmz
=\frac12+\int_0^{-vt^{-1/4}}g_1(r)\,dr.
\end{align}
We state several properties of the integration of $g_1(z)=\phi_1(|z|)$ based on Ferrero--Gazzola--Grunau\,\cite{fer;gaz;gru;08} 
and Gazzola--Grunau\,\cite{gaz;gru;09}. For $n\in\mathbb{Z}_+$, set
\[
\Phi_{N,n}(r):=\sum_{\ell=0}^n(-1)^\ell b_{N,\ell}\,r^{2\ell}, \quad
b_{N,\ell}:=\frac1{2^{N+1}\pi^{N/2}}\cdot\frac{\Gamma(\ell/2+N/4)}{2^{2\ell}\Gamma(\ell+1)\Gamma(\ell+N/2)}.
\]
The following lemma holds. 

\begin{lem}[{\cite[Lemma A.2]{fer;gaz;gru;08}}]\label{lem:fggA2}
Set $\delta_{N,n}:=b_{N,n+1}/b_{N,n}$ and assume that $n\in\mathbb{Z}_+$ is even. Then, for 
$0\le r\le1/\sqrt{\,\delta_{N,n}}$
\begin{align*}
&\Phi_{N,n-1}(r)\le\phi_N(r)\le\Phi_{N,n}, \\
&\max\{|\phi_N(r)-\Phi_{N,n}(r)|,|\phi_N(r)-\Phi_{N,n-1}(r)|,|\Phi_{N,n-1}(r)-\Phi_{N,n}(r)|\}\le b_{N,n}\,r^{2n}.
\end{align*}
\end{lem}

Set
\[
\Psi(r):=\int_0^r\phi_1(\eta)\,d\eta.
\]
Since $\phi_1$ is represented as a power series and converges locally uniformly in $\mathbb{R}$, we readily see that 
\begin{equation}\label{Psi_series}
\Psi(r)=\sum_{\ell=0}^\infty(-1)^\ell b_{1,\ell}\int_0^r\eta^{2\ell}\,d\eta
=\sum_{\ell=0}^\infty\frac{(-1)^\ell}{2\ell+1}b_{1,\ell}\,r^{2\ell+1}.
\end{equation}
For $n\in\mathbb{Z}_+$, define $\Psi_n(r)$ as 
\[
\Psi_n(r):=\int_0^r\Phi_{1,n}(\eta)\,d\eta=\sum_{\ell=0}^n\frac{(-1)^\ell}{2\ell+1}b_{1,\ell}\,r^{2\ell+1}.
\]
As a direct consequence of Lemma \ref{lem:fggA2}, we obtain the following corollary. 

\begin{cor}\label{cor:fggA2}
Assume that $n\in\mathbb{Z}_+$ is even. Then, for $0\le r\le1/\sqrt{\,\delta_{1,n}}$, 
\begin{align*}
&\Psi_{n-1}(r)\le\int_0^r\phi_1(\eta)\,d\eta\le \Psi_n(r), \\
&\max\biggl\{\biggl|\int_0^r\phi_1(\eta)\,d\eta-\Psi_n(r)\biggr|,\biggl|\int_0^r\phi_1(\eta)\,d\eta-\Psi_{n-1}(r)\biggr|,
\biggl|\Psi_n(r)-\Psi_{n-1}(r)\biggr|\biggr\} \\
&\hspace*{275pt}\le\frac{b_{1,n}}{2n+1}\,r^{2n+1}.
\end{align*}
\end{cor}

According to \cite[Theorem 2.3]{fer;gaz;gru;08}, $\phi_N$ changes its sign infinitely many times. Let 
$\{r_k^\pm\}_{k\in\mathbb{N}}$ be a sequence satisfying 
\begin{align*}
&\phi_1(r_k^\pm)=0, \quad 0<r_1^+<r_1^-<r_2^+<r_2^-<\cdots<r_k^+<r_k^-<\cdots, \\
&\phi_1(r)\left\{\begin{array}{ll} >0&\text{for}\,\ r\in\bigcup\limits_{k\in\mathbb{Z}_+}(r_k^-,r_{k+1}^+), \\
<0&\text{for}\,\ r\in\bigcup\limits_{k\in\mathbb{N}}(r_k^+,r_k^-), \end{array}\right.
\end{align*}
where $r_0^-=0$. Applying Lemma \ref{lem:fggA2} with $n=16$, we obtain
\begin{equation}\label{est_zeropt}
3.453<r_1^+<3.454, \quad 6.784<r_1^-<6.785.
\end{equation}
The following lemma can be proved by using \cite[Theorem 1 and Remark 1]{gaz;gru;09}. 

\begin{lem}\label{lem:Psi_properties}
\begin{list}{}{\leftmargin=0.9cm\labelwidth=1cm}
\item[(i)\,\,]
$\Psi(r)>0$ for $r>0$. 
\item[(ii)\,]
$\Psi(r)$ takes local maximum (resp. local minimum) at $r_k^+$\,(resp. $r_k^-$). 
\item[(iii)]
$\Psi(r)$ is strictly increasing (resp. decreasing) in each interval $(r_k^-,r_{k+1}^+)$ for $k\in\mathbb{Z}_+$ 
(resp. $(r_k^+,r_k^-)$ for $k\in\mathbb{N}$), where $r_0^-=0$. 
\end{list}
\end{lem}

\noindent
Furthermore, we are able to prove the following lemma. 

\begin{lem}\label{lem:monotonicity}
$\{\Psi(r_k^+)\}_{k\in\mathbb{N}}$ (resp. $\{\Psi(r_k^-)\}_{k\in\mathbb{N}}$) is strictly decreasing (resp. increasing). 
\end{lem}
\noindent
{\bf Proof.} 
It follows from \cite[Theorem 1]{gaz;gru;09} that for each $k\in\mathbb{N}$ 
\begin{equation}\label{intphi_sign}
\int_{r_k^-}^{r_{k+1}^-}\phi_1(r)\,dr>0, \quad \int_{r_k^+}^{r_{k+1}^+}\phi_1(r)\,dr<0.
\end{equation}
This implies
\begin{align*}
&\Psi(r_{k+1}^-)=\Psi(r_k^-)+\int_{r_k^-}^{r_{k+1}^-}\phi_1(r)\,dr>\Psi(r_k^-), \\
&\Psi(r_{k+1}^+)=\Psi(r_k^+)+\int_{r_k^+}^{r_{k+1}^+}\phi_1(r)\,dr<\Psi(r_k^+)
\end{align*}
for each $k\in\mathbb{N}$. \qed

\bigskip\noindent
With the help of \eqref{est_zeropt}, Corollary \ref{cor:fggA2} and some numerical computations,  
we observe that 
\begin{equation}\label{est_intphi}
\Psi(r_1^+/3)>0.32584, \quad 0.5522<\Psi(r_1^+)<0.5523, \quad 0.4938<\Psi(r_1^-)<0.4939
\end{equation}
For $ r\geq 0 $, set 
\[
I(r):=\Psi(r/(3a))-3\Psi(r/(2a))+3\Psi(r/a).
\]
Using \eqref{est_intphi}, we can prove the following lemma. 

\begin{lem}\label{lem:sign_I}
$I(r)>0$ for $r>0$. 
\end{lem}
\noindent
{\bf Proof.} 
We divide the proof into several cases. 

{\it Case 1}: $0<r/a\le r_1^+$. 
We readily see that 
\[
I(r)=\left(\int_0^{r/(3a)}+3\int_{r/(2a)}^{r/a}\right)\phi_1(\eta)\,d\eta>0.
\]

{\it Case 2}: $r_1^+<r/a\le r_1^-$. 
Note that $r_1^+/3<r/(3a)<r/(2a)\le r_1^-/2<r_1^+$ by \eqref{est_zeropt}. Then it follows from Lemma 
\ref{lem:Psi_properties} (iii) and \eqref{est_intphi} that 
\[
I(r)>\Psi(r_1^+/3)-3\Psi(r_1^+)+3\Psi(r_1^-)>0.32584+3(-0.5523+0.4938)=0.1329>0.
\]

{\it Case 3}: $r_1^-<r/a\le3r_1^-$. 
First, we derive the lower bound of $\Psi(r/(3a))$. 
Note that $r_1^+/3<r_1^-/3<r/(3a)\le r_1^-$. If $r/(3a)<r_1^+$, we have
\[
\Psi(r/(3a))=\Psi(r_1^+/3)+\int_{r_1^+/3}^{r/(3a)}\phi_1(\eta)\,d\eta>\Psi(r_1^+/3).
\]
If $r/(3a)\ge r_1^+$, Lemma \ref{lem:Psi_properties} (iii)  implies that $\Psi(r/(3a))\ge\Psi(r_1^-)$. 
Second, let us consider the upper bound of $\Psi(r/(2a))$. If $r/(2a)\le r_1^+$, it follows from Lemma 
\ref{lem:Psi_properties} (iii) 
that $\Psi(r/(2a))\le\Psi(r_1^+)$. If $r/(2a)>r_1^+$, we can choose $k_\ast\in\mathbb{N}$ such that 
$r_{k_\ast}^+<r/(2a)\le r_{k_\ast+1}^+$. 
Then, by \eqref{intphi_sign}, we see that 
\[
\Psi(r/(2a))=\Psi(r_1^+)+\sum_{k=1}^{k_\ast-1}\int_{r_k^+}^{r_{k+1}^+}\phi_1(\eta)\,d\eta
+\int_{r_{k_\ast^+}}^{r/(2a)}\phi_1(\eta)\,d\eta<\Psi(r_1^+).
\]
Finally, we derive the lower bound of $\Psi(r/a)$.  Choose $\ell_\ast\in\mathbb{N}$ such that 
$r_{\ell_\ast}^-<r/a\le r_{\ell_\ast+1}^+$. 
By virtue of \eqref{intphi_sign}, we obtain
\[
\Psi(r/a)=\Psi(r_1^-)+\sum_{\ell=1}^{\ell_\ast-1}\int_{r_\ell^-}^{r_{\ell+1}^-}\phi_1(\eta)\,d\eta
+\int_{r_{\ell_\ast}^-}^{r/a}\phi_1(\eta)\,d\eta
>\Psi(r_1^-). 
\]
Consequently, based on these bounds and \eqref{est_intphi}, we obtain 
\begin{align*}
I(r)>&\,\min\bigl\{\Psi(r_1^+/3),\Psi(r_1^-)\bigr\}-3\Psi(r_1^+)+3\Psi(r_1^-) \\
>&\,0.32584+3(-0.5523+0.4938)=0.1329>0.
\end{align*}

{\it Case 4}: $3r_1^-<r/a$. 
Choose $k_\ast\in\mathbb{N}$ such that $r_{k_\ast}^-\le r/(3a)<r_{k_{\ast}+1}^-$. Then we have
\[
\Psi(r/(3a))\left\{\begin{array}{ll} >\Psi(r_{k_\ast}^-)&\text{if}\,\ r/(3a)\in(r_{k_\ast}^-,r_{k_{\ast}+1}^+), \\ 
\ge\Psi(r_{k_\ast+1}^-)&\text{if}\,\ r/(3a)\in[r_{k_\ast+1}^+,r_{k_{\ast}+1}^-). \end{array}\right.
\]
Since $\Psi(r_k^-)\ge\Psi(r_1^-)$ for $k\in\mathbb{N}$ by Lemma \ref{lem:monotonicity}, we see that 
$\Psi(r/(3a))>\Psi(r_1^-)$. 
To derive the lower bound of $-\Psi(r/(2a))+\Psi(r/a)$, take $\ell_\ast,m_\ast\in\mathbb{N}$ satisfying 
\[
r_{\ell_\ast}^+\le r/(2a)<r_{\ell_\ast+1}^+, \quad r_{m_\ast}^-\le r/a<r_{m_\ast+1}^-.
\]
Then it follows from Lemma \ref{lem:monotonicity} that 
\[
-\Psi(r/(2a))+\Psi(r/a)>-\Psi(r_{\ell_\ast}^+)+\Psi(r_{m_\ast}^-)\ge-\Psi(r_1^+)+\Psi(r_1^-).
\]
Therefore, by virtue of these bounds and \eqref{est_intphi}, we see that 
\[
I(r)>\Psi(r_1^-)+3\{-\Psi(r_1^+)+\Psi(r_1^-)\}>0.4938+3(-0.5523+0.4938)=0.3183>0.
\]
This completes the proof. \qed

\medskip\noindent
Now we are ready to prove the Proposition \ref{prop:refinement}. 

\medskip\noindent
{\bf Proof of Proposition \ref{prop:refinement}.} 
For $\bmy\in\partial\Omega(t)$, take $\bmx\in\partial\Omega_0$ such that 
$|\bmy-\bmx|=|d(\bmy,\partial\Omega_0)|$. 
Set $v=d(\bmy,\partial\Omega_0)$. Then we have $\bmy=\bmx+v\bmn(\bmx)$. 
Applying Theorem \ref{thm:asympt-exp_4w2}, we obtain 
\begin{align*}
&\hspace*{-10pt}
u(\bmx+v\,\bmn(\bmx),a^4t) \\
=&\,\frac12+\sum_{m=0}^\infty\frac{a^{2m}(-\lambda)^mt^{m/2}}{m!}
\int_{\mathbb{R}^{N-1}}\int_0^{\psi_a(\sbmz',v,t)}(-\Delta_{\sbmz})^mg_N(\bmz)\,dz_Nd\bmz'
+O_{\p\om_0}(e^{-\mu_\ast t^{-1/3}}),
\end{align*}
where 
\begin{equation}\label{def_psia}
\psi_a(\bmz',v,t):=\frac1{at^{1/4}}\{-v+f(at^{1/4}\bmz')\}
\end{equation}
for $a>0$ and a function $f$ satisfying (A1)--(A3). Then we have  
\begin{align*}
&\hspace*{-10pt}
U(\bmx+v\,\bmn(\bmx),t) \\
=&\,\frac12+\sum_{m=0}^\infty\frac{a^{2m}(-\lambda)^mt^{m/2}}{m!} \\
&\hspace*{40pt}
   \cdot\left(3^{2m}\int_{\mathbb{R}^{N-1}}\int_0^{\psi_{3a}(\sbmz',v,t)}-3\cdot2^{2m}\int_{\mathbb{R}^{N-1}}
\int_0^{\psi_{2a}(\sbmz',v,t)}
+3\int_{\mathbb{R}^{N-1}}\int_0^{\psi_a(\sbmz',v,t)}\right) \\
&\hspace*{285pt}
(-\Delta_{\sbmz})^mg_N(\bmz)\,dz_Nd\bmz' \\
&\,+O_{\p\om_0}(e^{-\mu_\ast t^{-1/3}}).
\end{align*}
%
%
Since $\bmy\in\partial\Omega(t)$ and $v=d(\bmy,\partial\Omega_0)$, we see that 
\[
U(\bmx+v\,\bmn(\bmx),t)=\frac12.
\]
Furthermore, Proposition \ref{prop:inclusion} implies that $|v|\le K_\ast t^{1/4}$ for $t\in(0,t_\ast)$. 

Let us prove $v=o_{\text{{\tiny$\partial\Omega_0$}}}(t^{1/4})$ as $t\to+0$. 
In the proof by contradiction, suppose that there exists $\varepsilon_\ast>0 $ 
such that for any $ k\in\bN $ there are $ \{t_k\}_{k\in\bN} $ and 
$ \{\bmy_k\}_{k\in\bN}\subset\partial\Omega(t_k) $  
such that 
\[
0<t_k<\frac1k, \quad |d(\bmy_k,\partial\Omega_0)|>\varepsilon_\ast t_k^{1/4}. 
\]
Then there exist sequences $\{\bmx_k\}_{k\in\mathbb{N}}\subset\partial\Omega_0$ and 
$\{r_k\}_{k\in\mathbb{N}}$ satisfying 
\[
\bmy_k=\bmx_k+v_k\bmn_k,\quad \varepsilon_\ast<|r_k|\le K_\ast, \quad U(\bmx_k+v_k\,\bmn_k,t_k)=\frac12,
\]
where $v_k:=-r_kt_k^{1/4}$ and $\bmn_k:=\bmn(\bmx_k)$. Since $\{\bmx_k\}_{k\in\mathbb{N}}$ and $\{r_k\}_{k\in\mathbb{N}}$ 
are bounded sequences, by taking subsequences if necessary, we may assume that 
$ \bmx_k\to\bmx_\ast\in\p\om_0 $, $r_k\to r_\ast$ as $k\to\infty$ where $\varepsilon_\ast\le|r_\ast|\le K_\ast$. 
On the other hand, $\psi_a(\bmz',v_k,t_k)$ is represented as 
\begin{align*}
\psi_a(\bmz',v_k,t_k)
=&\,(at_k^{1/4})^{-1}\left\{
-v_k+\frac{(at_k^{1/4})^2\langle\bmz',\nabla_{\sbmz'}\rangle_{N-1}^2f(\theta at_k^{1/4}\bmz')}2\right\} \\
=&\,\frac{r_k}a+\frac{at_k^{1/4}}2\langle\bmz',\nabla_{\sbmz'}\rangle_{N-1}^2f(\theta at_k^{1/4}\bmz')
\end{align*}
for some $\theta\in(0,1)$. Since $\psi_a(\bmz',v_k,t_k)\to r_\ast/a$ as $k\to\infty$, 
it follows from \eqref{sol_wholesp} that 
\begin{align*}
\frac12
&=\lim_{k\to\infty}U(\bmx_k+v_k\,\bmn_k,t_k) \\ \nonumber
&=\frac12+
\left(\int_{\mathbb{R}^{N-1}}\int_0^{r_*/3a}-3\int_{\mathbb{R}^{N-1}}
\int_0^{r_*/2a}+3\int_{\mathbb{R}^{N-1}}\int_0^{r_*/a}\right)g_N(\bmx)\,d\bmx \\ \nonumber
&=\frac12+\left(\int_0^{r_\ast/(3a)}-3\int_0^{r_\ast/(2a)}+3\int_0^{r_\ast/a}\right)g_1(r)\,dr.
\end{align*}
%
By Lemma \ref{lem:sign_I} and $ g_1(r)=\varphi_1(|r|)$, we see that if $r_\ast>0$
\[
\left(\int_0^{r_\ast/(3a)}-3\int_0^{r_\ast/(2a)}+3\int_0^{r_\ast/a}\right)g_1(r)\,dr>0,
\]
and if $r_\ast<0$
\begin{align*}
&\left(\int_0^{r_\ast/(3a)}-3\int_0^{r_\ast/(2a)}+3\int_0^{r_\ast/a}\right)g_1(r)\,dr \\
&=-\left(\int_0^{-r_\ast/(3a)}-3\int_0^{-r_\ast/(2a)}+3\int_0^{-r_\ast/a}\right)g_1(r)\,dr<0.
\end{align*}
These facts lead to a contradiction. Therefore, $v=o_{\partial\Omega_0}(t^{1/4})$ as $t\to+0$. \qed

\bigskip
We next prove an estimate of the derivative of $U$, which guarantees that $\partial\Omega(t)$ is a smooth hypersurface. 
Set $M_{0,a}:=11M_0/(6a)$, where $M_0$ is as in Lemma \ref{lem:moment}.

\begin{prop}\label{prop:derivative}
There exist $\varepsilon_0\in(0,1)$ and $C>0$ such that for $t>0$ small enough, $\bmx\in\partial\Omega_0$ and 
$v\in(-\delta_0,\delta_0)$ satisfying $|v|\le\varepsilon_0t^{1/4}$ 
\begin{align}\label{est-derivative}
&-\frac{3M_{0,a}}2\,t^{-1/4}-C(\varepsilon_0+t^{1/4}) \\
&\le\langle\nabla_{\sbmx}U(\bmx+v\,\bmn(\bmx),t),\bmn(\bmx)\rangle_N\le-\frac{M_{0,a}}2\,t^{-1/4}
+C(\varepsilon_0+t^{1/4}). \nonumber
\end{align}
Here $ \delta_0\in (0,1/2) $ has been taken at the end of section 3.1. 
Furthermore, $\partial\Omega(t)$ is a smooth hypersurface for $t>0$ small enough. 
\end{prop}
\noindent
{\bf Proof.} 
{\it Step 1}. Set 
\[
\Xi_{N,a}^{(1)}:=-\frac1{at^{1/4}}\int_{\mathbb{R}^{N-1}}\int_{-\infty}^{\psi_a(\sbmz',v,t)}\partial_{z_N}
g_N(\bmz)\,dz_Nd\bmz'.
\]
Then we have 
\begin{equation}\label{nabla_u_normal_1}
|\langle\nabla_{\sbmx}u_a(\bmx+v\,\bmn(\bmx),t),\bmn(\bmx)\rangle_N-\Xi_{N,a}^{(1)}|\le C_1\,|\lambda|\,t^{1/4}.  
\end{equation}
Indeed, applying an argument similar to the proof of Lemma \ref{lem:series_finite}, we obtain 
%
%
\[
\left|\frac1{at^{1/4}}\sum_{m=1}^\infty\frac{a^{2m}(-\lambda)^mt^{m/2}}{m!}
\int_{\mathbb{R}^{N-1}}\int_{-\infty}^{\psi_a(\sbmz',v,t)}\partial_{z_N}(-\Delta_{\sbmz})^mg_N(\bmz)\,dz_Nd\bmz'\right|
\le C_1\,|\lambda|\,t^{1/4}
\]
for $|\lambda|t^{1/2}\le\gamma/a^2$, where $\psi_a(\bmz',v,t)$ is given by \eqref{def_psia}. 
This inequality and \eqref{Dsol_expansion} yield \eqref{nabla_u_normal_1}. 

{\it Step 2}. We prove (\ref{est-derivative}). 
By \eqref{threshold-func} and \eqref{nabla_u_normal_1}, there exists a $C_2>0$ such that 
\begin{equation}\label{nabla_U_normal}
|\langle\nabla_{\sbmx}U(\bmx+v\,\bmn(\bmx),t),\bmn(\bmx)\rangle_N-(\Xi_{N,3a}^{(1)}-3\Xi_{N,2a}^{(1)}+3\Xi_{N,a}^{(1)})|
\le C_2\,|\lambda|\,t^{1/4}
\end{equation}
for small $t>0$.
Note that $\Xi_{N,a}^{(1)}$ is rewritten as
\[
\Xi_{N,a}^{(1)}=-\frac1{at^{1/4}}\int_{\mathbb{R}^{N-1}}g_N(\bmz',\psi_a(\bmz',v,t))\,d\bmz'.
\]
It follows from \eqref{Dg_zN=0} and the Taylor's theorem that 
\[
g_N(\bmz',\psi_a(\bmz',v,t))
=g_N(\bmz',0)+\frac12\partial_{z_N}^2g_N(\bmz',\theta\psi_a(\bmz',v,t))\{\psi_a(\bmz',v,t)\}^2
\]
for some $\theta\in(0,1)$. This and the definition of $M_0$ imply that
\begin{equation}\label{Xi_Na}
\Xi_{N,a}^{(1)}
=-\frac{M_0}{at^{1/4}}
-\frac1{2at^{1/4}}\int_{\mathbb{R}^{N-1}}\partial_{z_N}^2g_N(\bmz',\theta\psi_a(\bmz',v,t))\{\psi_a(\bmz',v,t)\}^2\,d\bmz'.
\end{equation}
In addition, taking account of 
\[
\psi_a(\bmz',v,t)=-\frac{v}{at^{1/4}}+\frac{at^{1/4}}2
\langle\bmz',\nabla_{\sbmz'}\rangle_{N-1}^2f(\widetilde{\theta}at^{1/4}\bmz')
\]
for some $\widetilde{\theta}\in(0,1)$, we have
\begin{align}
\label{eq040201}
&\hspace*{-10pt}
\int_{\mathbb{R}^{N-1}}\partial_{z_N}^2g_N(\bmz',\theta\psi_a(\bmz',v,t))\{\psi_a(\bmz',v,t)\}^2\,d\bmz' \\
\nonumber
=&\,\frac{v^2}{(at^{1/4})^2}\int_{\mathbb{R}^{N-1}}\partial_{z_N}^2g_N(\bmz',\theta\psi_a(\bmz',v,t))\,d\bmz'\\
\nonumber
&\,-v\int_{\mathbb{R}^{N-1}}\partial_{z_N}^2g_N(\bmz',\theta\psi_a(\bmz',v,t))
\langle\bmz',\nabla_{\sbmz'}\rangle_{N-1}^2f(\widetilde{\theta}at_k^{1/4}\bmz')\,d\bmz' \\
\nonumber
&\,+\frac{(at^{1/4})^2}4\int_{\mathbb{R}^{N-1}}\partial_{z_N}^2g_N(\bmz',\theta\psi_a(\bmz',v,t))
\bigl\{\langle\bmz',\nabla_{\sbmz'}\rangle_{N-1}^2f(\widetilde{\theta}at_k^{1/4}\bmz')\bigr\}^2\,d\bmz'.
\end{align}
We assume that $|v|\le\varepsilon_0t^{1/4}$ where 
$\varepsilon_0\in(0,1)$ is suitably chosen later. Then it follows that 
\[
|\psi_a(\bmz',v,t)|\le C_3(|\bmz'|^2t^{1/4}+\varepsilon_0).
\]
This fact and Corollary \ref{cor:Cui-thm} with $ |\al|=2 $ and $m=0$ yield that 
\begin{align*}
|\partial_{z_N}^2g_N(\bmz',\theta\psi_a(\bmz',v,t))|
\le&\,C_4\bigl\{1+(|\bmz'|^2+|\psi_a(\bmz',v,t)|^2)^{1/2}\bigr\}^{2/3}e^{-\mu|\sbmz'|^{4/3}} \\
\le&\,C_4(1+|\bmz'|+|\psi_a(\bmz',v,t)|)^{2/3}e^{-\mu|\sbmz'|^{4/3}} \\
\le&\,C_5(1+|\bmz'|^{2/3}+|\bmz'|^{4/3}t^{1/6})e^{-\mu|\sbmz'|^{4/3}}.
\end{align*}
Applying this estimate to the right-hand side of (\ref{eq040201}), we are able to find a constant 
$C_6>0$ such that 
\begin{align*}
&\hspace*{-10pt}
\left|\frac12\int_{\mathbb{R}^{N-1}}\partial_{z_N}^2
g_N(\bmz',\theta\psi_a(\bmz',v,t))\{\psi_a(\bmz',v,t)\}^2\,d\bmz'\right|
\le C_6\left\{\frac{\varepsilon_0^2}{a^2}+\varepsilon_0t^{1/4}+(at^{1/4})^2\right\}
\end{align*}
for $t>0$ small enough. Recalling \eqref{Xi_Na}, we see that
\[
\left|\,\Xi_{N,3a}^{(1)}-3\Xi_{N,2a}^{(1)}+3\Xi_{N,a}^{(1)}+M_{0,a}\,t^{-1/4}\right|
\le C_7\left(\frac{\varepsilon_0^2}{a^3t^{1/4}}+\frac{\varepsilon_0}{a}+at^{1/4}\right)
\]
for some $ C_7>0 $. Choosing $\varepsilon_0>0$ such that
\[
\frac{C_7\varepsilon_0^2}{a^2}\le\frac{M_0}2,
\]
we obtain
\begin{align*}
&-\frac{3M_{0,a}}{2}t^{-1/4}-C_7\left(\frac{\varepsilon_0}{a}+at^{1/4}\right) \\
&\le \Xi_{N,3a}^{(1)}-3\Xi_{N,2a}^{(1)}+3\Xi_{N,a}^{(1)} 
\le-\frac{M_{0,a}}{2}t^{-1/4}+C_7\left(\frac{\varepsilon_0}{a}+at^{1/4}\right).
\end{align*}
By this inequality and \eqref{nabla_U_normal}, we are led to the desired result. 

{\it Step 3}. Let us prove that $\partial\Omega(t)$ is a smooth hypersurface for $t>0$ small enough. 
By virtue of \eqref{est-derivative}, $\langle\nabla_{\sbmx}U(\bmx+v\,\bmn(\bmx),t),\bmn(\bmx)\rangle_N$ is far from zero 
for $t>0$ small enough, in particular, it is negative. This fact and the implicit function theorem implies the desired result. 
\qed


\begin{thm}\label{thm:est-distance}
There exists a $C_0>0$ such that  
\[
     \sup_{\sbmx\in\partial\Omega(t)}|d(\bmx,\partial\Omega_0)|\le C_0\,t
\]
for $t>0$ small enough. 
\end{thm}
\noindent
{\bf Proof.} 
Applying \eqref{asympt_expansion} with $V=0$, there exists a $C_1>0$ such that 
\begin{equation}\label{est_U-1/2}
\left|U(\bmx,t)-\frac12\right|\le C_1t^{3/4}
\end{equation}
for $\bmx\in\partial\Omega_0$ and $t>0$ small enough. For any $\bmy\in\partial\Omega(t)$, let 
$\bmx\in\partial\Omega_0$ 
be a point satisfying $|d(\bmy,\partial\Omega_0)|=|\bmy-\bmx|$. Then $\bmy$ can be represented as 
$\bmy=\bmx+v\bmn(\bmx)$ 
where $v=d(\bmy,\partial\Omega_0)$. This implies that
\[
U(\bmx+v\bmn(\bmx),t)-U(\bmx,t)=\langle\nabla_{\sbmx} U(\bmx+\theta v\bmn(\bmx),t),\bmn(\bmx)\rangle_N\,v
\]
for some $\theta\in(0,1)$. Taking account of $U(\bmy,t)=1/2$ for $\bmy\in\partial\Omega(t)$ and using 
Proposition \ref{prop:derivative} 
and \eqref{est_U-1/2}, we see that 
\begin{align*}
&|\langle\nabla_{\sbmx} U(\bmx+\theta v\bmn(\bmx),t),\bmn(\bmx)\rangle_N\,v|\le C_1\,t^{3/4}, \\
&-C_2\,t^{-1/4}\le\langle\nabla_{\sbmx} U(\bmx+\theta v\bmn(\bmx),t),\bmn(\bmx)\rangle_N\le -C_3\,t^{-1/4}
\end{align*}
for $t>0$ small enough where $C_2,C_3$ are positive constants independent of  
$ \bmx\in\partial\Omega_0 $ and $t$. As a result, it follows that 
there exists $C_4>0$ such that $|v|\le C_4\,t$ for $ \bmx\in\partial\Omega_0 $ and small 
$t>0$. This is the desired result. \qed

\bigskip\noindent
From this theorem, it follows that $ V $, defined as in \eqref{eq:def-V}, is bounded on $ \p\om_0\times(0,t_0) $ for some small 
$ t_0>0 $ and hence the argument in subsection \ref{subsec:proposed_algorithm} is justified, that is:

\begin{thm}\label{th040201}
Let $V$ be as in \eqref{eq:def-V}. Then there exist $ C>0 $ and $ t_0>0 $ such that
\[
\left|\,V+\Delta_gH+H|A|^2-\frac12H^3
+2\sum_{\substack{i_1,i_2,i_3\in\Lambda \\ i_1<i_2<i_3}}\kappa_{i_1}\kappa_{i_2}\kappa_{i_3}-\lambda H\,\right|
\leq Ct^{1/4}
\]
for all $ t\in (0,t_0) $ and $ \bmx\in\p\om_0 $. 
Especially, this estimate turns to 
\[
\begin{array}{ll}
\displaystyle\left|\,V+\kappa_{ss}+\frac12\kappa^3-\lambda \kappa\,\right|\leq Ct^{1/4}&\text{if}\,\ N=2, \\[0.45cm]
\displaystyle\left|\,V+\Delta_gH+H|A|^2-\frac12H^3-\lambda H\,\right|\leq Ct^{1/4}&\text{if}\,\ N=3
\end{array}
\]
for all $ t\in (0,t_0) $ and $ \bmx\in\p\om_0 $. 
\end{thm}


\begin{rem}
(i) If $ \p\om_0 $ is of class $ C^n\,(n\geq 6)$, we can replace $ t^{1/4} $ with $ t^{1/2} $ in 
Theorem \ref{th040201}.  
But we omit the detail because the actual calculations are more complicated.   

(ii) Based on the boundedness of $ V $, we can improve (\ref{est-derivative}) as follows:
$$
   |\langle\nabla_{\sbmx}U(\bmx+v\,\bmn(\bmx),t),\bmn(\bmx)\rangle_N+M_{0,a}t^{-1/4}|
   \leq Ct^{3/4}.  
$$
Moreover, we can also estimate 
$$
   |\langle\nabla_{\sbmx}U(\bmx+v\,\bmn(\bmx),t),\bmtau(\bmx)\rangle|\leq Ct^{3/4}   
$$
where $ \bmtau(\bmx) $ is any unit tangential vector of $ \p\om_0 $ at $ \bmx\in\p\om_0 $.  
Thus the outer unit normal $ -\nabla_{\sbmx}U(\cdot,t)\,/\,|\nabla_{\sbmx} U(\cdot,t)| $ 
of $ \p\om(t) $ is nearly equal to 
$ \bmn(\cdot) $, that of $ \p\om_0 $, for any small $ t>0 $.  

\end{rem}
\section{Numerical experiments}\label{sec:numerical_experiments}
In this section, we present some results of numerical experiments based on our thresholding algorithm 
which is proposed in subsection \ref{subsec:proposed_algorithm}.
We focus on the case where $N=2$ and $ \lambda=0 $, that is, the Willmore flow for planar curves.
\subsection{Numerical scheme}
Let $D=(\mathbb{R}/L\mathbb{Z})\times(\mathbb{R}/L\mathbb{Z})$ be a periodic square region with $L>0$ and 
assume that $\Omega_0\subset D$ holds. Then we solve the following initial value problem
\begin{equation}\label{eq:target}
    \left\{\begin{array}{l}
    u_t=-(-\Delta)^2u\,\,\ \text{in}\ D\times(0,\infty), \\
    u(\cdot,0)=\chi_{\Omega_0}\,\,\ \text{in}\ D.
    \end{array}\right.
\end{equation}
By the definition of $D$, the periodic boundary condition is implicitly imposed on the above problem.
Note that there are other possibilities regarding the boundary conditions, for example, the Dirichlet boundary condition 
or the Neumann boundary condition is natural.
However, we adopt the periodic boundary condition here since we develop a numerical scheme based on the Fourier transform.

We seek a solution in the form of Fourier series
\begin{equation*}
    u(\bmx,t)=\sum_{\sbmxi\in\mathbb{Z}^2}u_{\sbmxi}(t)e^{2\pi\textbf{i}\langle\sbmx,\sbmxi\rangle_2/L},
\end{equation*}
since $u$ is supposed to be a periodic function.
Taking account of 
\begin{equation*}
    u_t=\sum_{\sbmxi\in\mathbb{Z}}\dot{u}_{\sbmxi}(t)e^{2\pi\textbf{i}\langle\sbmx,\sbmxi\rangle_2/L},
    \quad
    (-\Delta)^2u=\sum_{\sbmxi\in\mathbb{Z}}\frac{16\pi^4}{L^4}|\bmxi|^4u_{\sbmxi}(t)e^{2\pi\textbf{i}\langle\sbmx,\sbmxi\rangle_2/L}
\end{equation*}
and substituting these expressions into the first equation in \eqref{eq:target}, we are led to an infinite system 
of ordinary differential equations
\begin{equation*}
    \dot{u}_{\sbmxi}(t)=-\frac{16\pi^4}{L^4}|\bmxi|^4u_{\sbmxi}(t),\quad\bmxi\in\mathbb{Z}^2,
\end{equation*}
where the dot symbol indicates the time derivative.
Each ordinary differential equation can be separately and explicitly solved as
\begin{equation*}
    u_{\sbmxi}(t)=\chi_{\Omega_0,\sbmxi}\,e^{-16\pi^4|\sbmxi|^4t/L^4},
\end{equation*}
where $\chi_{\Omega_0,\sbmxi}$ denotes the Fourier coefficient of the characteristic function $\chi_{\Omega_0}$.
Then the function $U$ defined by \eqref{threshold-func} can be computed as
\begin{align}    \label{eq:U_Fourier}
    &U(\bmx,h)
    =u_{3a}(\bmx,h)-3u_{2a}(\bmx,h)+3u_a(\bmx,h) \\
    &=\sum_{\sbmxi\in\mathbb{Z}^2}\chi_{\Omega_0,\sbmxi}\left(
        e^{-16\pi^4|\sbmxi|^4(3a)^4h/L^4}-3e^{-16\pi^4|\sbmxi|^4(2a)^4h/L^4}+3e^{-16\pi^4|\sbmxi|^4a^4h/L^4}
    \right)e^{2\pi\textbf{i}\langle\sbmx,\sbmxi\rangle_2/L}, \notag
\end{align}
where $h>0$ denotes a time step.
Hence we obtain $\Omega_1$ by $\Omega(h)$, where $\Omega(t)$ is given by \eqref{set_Omega_t}.
Repeating the above procedure yields an approximation of the Willmore flow.

In the actual computation, we first generate a uniform mesh $\{\bmx_{ij}\}_{(i,j)\in\Delta_N^2}$ as
\begin{equation*}
    \bmx_{ij}=\begin{pmatrix}
        \left(i+\frac{1}{2}\right)r\\[1ex]
        \left(j+\frac{1}{2}\right)r
    \end{pmatrix},
    \quad
    (i,j)\in\Delta_N^2,
\end{equation*}
where $r:=L/N$, $\Delta_N:=\{0,1,\ldots,N-1\}$, and $N$ is a power of $2$.
Let us denote an approximate value of $u(\cdot,t)$ at $(x_i,y_j)$ by $u_{ij}(t)$.
Using the fast Fourier transform, we can compute discrete Fourier coefficients $\chi_{\Omega_0,\sbmxi}$ 
for $\bmxi\in\Delta_N^2$ with $O(N^2\log N)$ complexity. Then we are able to approximate 
$\{U(\bmx_{ij},h)\}_{(i,j)\in\Delta_N^2}$ given by \eqref{eq:U_Fourier} 
by means of applying the inverse fast Fourier transform to 
$\{\chi_{\Omega_0,\sbmxi}\left(u_{\sbmxi}(81a^4h)-3u_{\sbmxi}(16a^4h)+3u_{\sbmxi}(a^4h)\right)\}_{\sbmxi\in\Delta_N^2}$. 

\begin{rem}
    As described above, we use a uniform mesh to perform the fast Fourier transform and its inverse.
    It is preferable to adopt some mesh refinement techniques to capture the evolution of the interface more precisely.
    Indeed, Ruuth\,\cite{ruu;98} developed a gradually adaptive mesh refinement technique and applied it to the computation 
    of the curve shortening and mean curvature flows based on the BMO algorithm.
    This technique has also been applied to the Willmore flow in Esedo\=glu-Ruuth-Tsai\,\cite{ese;ruu;tsa;08}.
    It could be expected that a similar approach would be necessary for our algorithm. However, as a first step of numerical computation 
    of the Willmore flow based on our proposed scheme, we adopt the uniform mesh in this paper.
\end{rem}

\subsection{Numerical results}

In this subsection, we show two numerical results of our algorithm.

\subsubsection{Self-similar solution}
\begin{figure}[htbp]
    \centering
    \includegraphics[width=.485\hsize]{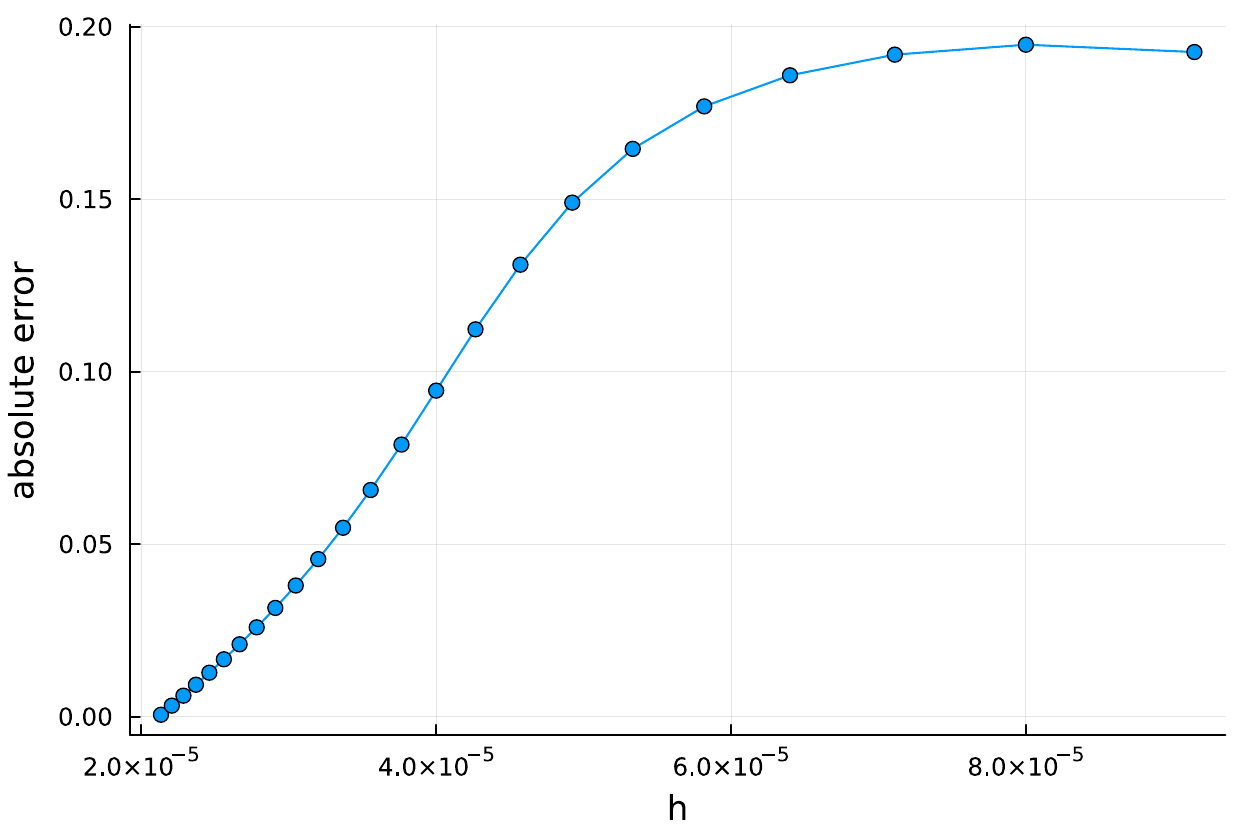}\,\,\ \includegraphics[width=.485\hsize]{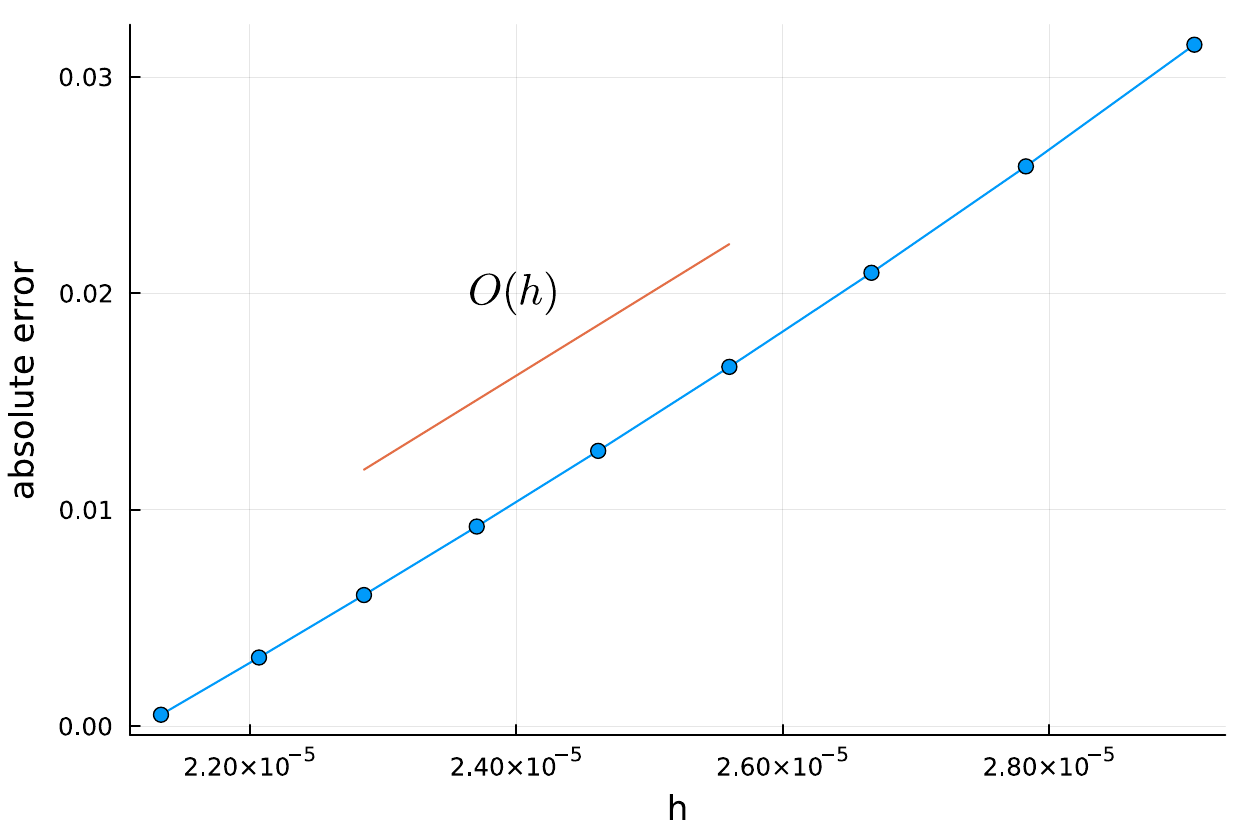}
    \caption{Plot of the absolute error versus the time step size $h$ at $t=0.00064$. 
    The left is the overview and the right is the close-up around $h=2.50\times10^{-5}$.}
    \label{fig:error}
\end{figure}

First, let us consider the case where the initial curve is a circle.
As is well documented, the solution develops over time, maintainig a circular shape.
Its radius, designated as $R$, is defined as the solution to the following ordinary differential equation:
\begin{equation*}
    \dot{R}=\frac{1}{2}\left(\frac{1}{R}\right)^3.
\end{equation*}
That is to say, $R$ is specifically represented as 
\begin{equation*}
    R(t)=\sqrt[4]{R(0)^4+2t}.
\end{equation*}
Figure~\ref{fig:error} illustrates how the discrepancy between the area of the exact solution and that of 
the numerical solution diminishes as $h$ is gradually reduced.
The initial radius is $0.1$, and the final computation time is $0.00064$.
The remaining parameters were determined as follows:
\begin{list}{$\bullet$}{\leftmargin=1cm\topsep=0.325cm\itemsep=0cm}
    \item $L=1$ (the domain size),
    \item $N=2^{14}$ (the number of mesh).
\end{list}
As illustrated in the graph, the approximation error exhibits a gradual decline as $h$ is reduced.
The orange line in the graph represents $O(h)$, indicating that the proposed scheme is first-order accurate.
This is a superior value compared to the result by Esedo\=glu--Ruuth--Tsai\,\cite{ese;ruu;tsa;08}, 
which suggested one-half order accuracy.

\subsubsection{Other nontrivial examples}

In this subsection, we show several nontrivial numerical results.
Parameters are chosen as follows:
\begin{list}{$\bullet$}{\leftmargin=1cm\topsep=0.325cm\itemsep=0cm}
    \item $L=5$ (the domain size),
    \item $N=2^{10}=1024$ (the number of mesh).
\end{list}
Note that the time step $h$ cannot simply be made smaller.
As pointed out by Ruuth\,\cite[Section 2.3]{ruu;98}, the spatial step must be sufficiently smaller than the time step 
to allow the interface to move.

\begin{figure}[htbp]
    \begin{center}
        \begin{minipage}{.3\hsize}
            \centering
            \includegraphics[width=\hsize, bb = 0 0 400 400]{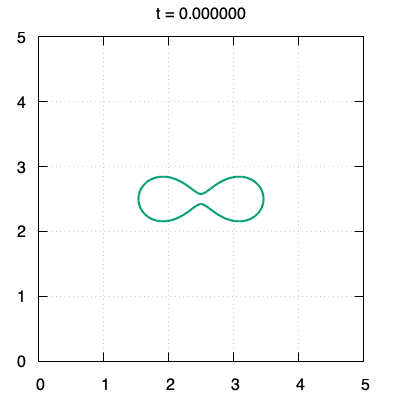}
        \end{minipage}%
        \begin{minipage}{.3\hsize}
            \centering
            \includegraphics[width=\hsize, bb = 0 0 400 400]{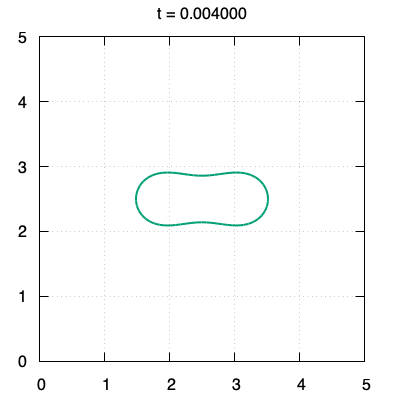}
        \end{minipage}%
    \\
        \begin{minipage}{.3\hsize}
            \centering
            \includegraphics[width=\hsize, bb = 0 0 400 400]{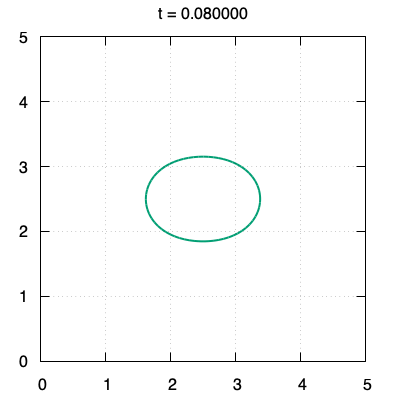}
        \end{minipage}%
        \begin{minipage}{.3\hsize}
            \centering
            \includegraphics[width=\hsize, bb = 0 0 400 400]{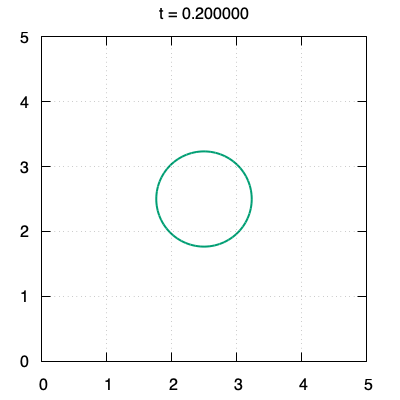}
        \end{minipage}%
        \caption{Numerical results of the Willmore flow with initial region given as the Cassini's 
        Parameters are chosen as $L=5, N=1024$, and $h=0.004$.}
        \label{fig:cassini}
    \end{center}
    \end{figure}

First, consider the case where the following Cassini's oval gives the initial region:
\begin{equation}
    \Omega_0:=\left\{
        \bmx\in\mathbb{R}^2
        \mid
        (x_1^2+x_2^2)^2-2b^2(x_1^2-x_2^2)\le a^4-b^4
    \right\},
    \,\,\ 
    a=0.6825,\ b=0.678.
    \label{eq:cassini}
\end{equation}
The results are shown in Figure~\ref{fig:cassini}, where the time step $h$ is chosen as $h=0.004$. 
It can be seen that the numerical computations are performed stably without any numerical instability.

\begin{figure}[htbp]
\begin{center}
    \begin{minipage}{.3\hsize}
        \centering
        \includegraphics[width=\hsize, bb = 0 0 400 400]{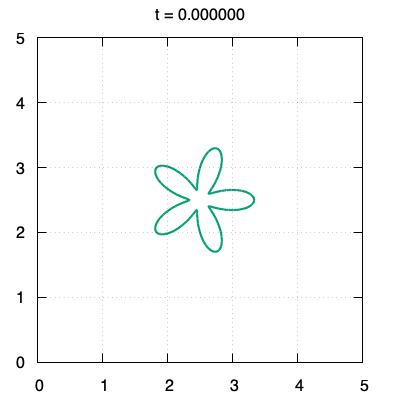}
    \end{minipage}%
    \begin{minipage}{.3\hsize}
        \centering
        \includegraphics[width=\hsize, bb = 0 0 400 400]{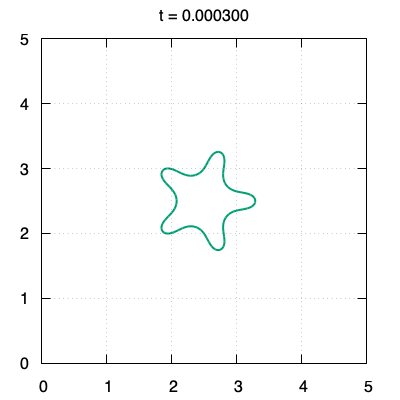}
    \end{minipage}%
\\
    \begin{minipage}{.3\hsize}
        \centering
        \includegraphics[width=\hsize, bb = 0 0 400 400]{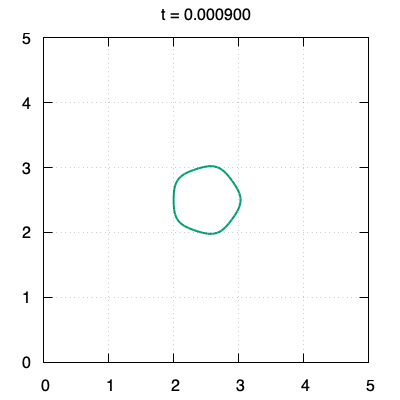}
    \end{minipage}%
    \begin{minipage}{.3\hsize}
        \centering
        \includegraphics[width=\hsize, bb = 0 0 400 400]{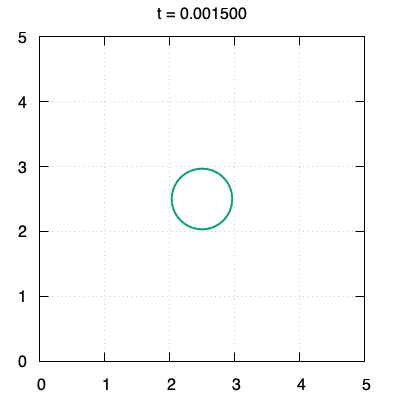}
    \end{minipage}%
    \caption{Numerical results of the Willmore flow with initial region given by \eqref{eq:rose}.
    Parameters are chosen as $L=5, N=1024$, and $h=0.0003$.}
    \label{fig:rose}
\end{center}
\end{figure}

Another more complex initial shape is the following initial region:
\begin{equation}
    \Omega_0:=\left\{
        \bmx\in\mathbb{R}^2
        \mid
        x_1^2+x_2^2\le\max\{0.01,r(\bmx)^2\}
    \right\},
    \label{eq:rose}
\end{equation}
where 
\begin{equation*}
    c=\frac{x_1}{\sqrt{x_1^2+x_2^2}},
    \quad
    r(\bmx)=0.5+\frac{16c^5-20c^3+5c}{3}.
\end{equation*}
The results are depicted in Figure~\ref{fig:rose}, in which the time step $h$ is chosen as $h=0.003$. 
In this case, the numerical computation is also successful without instability.

\appendix

\section{Mean curvature and its derivatives} \label{sec:mc}
Set $\Phi(\bmx')=(\bmx',f(\bmx'))$ where $\bmx'=(x_1,\cdots,x_{N-1})$. In this subsection, we denote $\partial_{x^i}$ 
by $\partial_i$. 
First, we have
\begin{align*}
&\partial_i\Phi=(0,\cdots,0,\underset{\hat{i}}{1},0,\cdots,0,\partial_if), \\
&\partial_1\Phi\times\cdots\times\partial_{N-1}\Phi=(-\nabla_{\sbmx^\prime} f,1).
\end{align*}
Then the outer unit normal $\bmn$ to the hypersurface $\{\Phi(\bmx')\,|\,\bmx'\in D\}$ for $D\subset\mathbb{R}^{N-1}$ 
is represented as
\[
\bmn=\frac{(-\nabla_{\sbmx^\prime} f,1)}{\sqrt{1+|\nabla_{\sbmx^\prime} f|^2}}.
\]
Also, the first fundamental forms are
\[
g_{ij}=\langle\partial_i\Phi,\partial_j\Phi\rangle_N=\delta_{ij}+\partial_if\,\partial_jf
\]
Setting $g=\det(g_{ij})=1+|\nabla_{\sbmx^\prime} f|^2$, we see that 
\[
g^{ij}=\delta_{ij}-\frac{\partial_if\,\partial_jf}{1+|\nabla_{\sbmx^\prime} f|^2}
\]
Let us derive the second fundamental forms. Since 
\[
\partial_i\partial_j\Phi=(0,\cdots,0,\partial_i\partial_jf),
\]
we are led to
\[
h_{ij}=\langle\partial_i\partial_j\Phi,\bmn\rangle_N=\frac{\partial_i\partial_jf}{\sqrt{1+|\nabla_{\sbmx^\prime} f|^2}}.
\]
Thus the mean curvature $H$ is represented as
\[
H=\sum_{i,j\in\Lambda}g^{ij}h_{ij}
=\sum_{i,j\in\Lambda}\left(\delta_{ij}-\frac{\partial_if\,\partial_jf}{1+|\nabla_{\sbmx^\prime} f|^2}\right)
\frac{\partial_i\partial_jf}{\sqrt{1+|\nabla_{\sbmx^\prime} f|^2}}, \\
\]
Furthermore, we have
\begin{align*}
|A|^2=&\,\sum_{i,j,k,\ell\in\Lambda}g^{ik}g^{j\ell}h_{ij}h_{k\ell} \\
=&\,\sum_{i,j,k,\ell\in\Lambda}\left(\delta_{ik}-\frac{\partial_if\,\partial_kf}{1+|\nabla_{\sbmx^\prime} f|^2}\right)
\left(\delta_{j\ell}-\frac{\partial_jf\,\partial_{\ell}f}{1+|\nabla_{\sbmx^\prime} f|^2}\right)
\frac{\partial_i\partial_jf\,\partial_k\partial_\ell f}{1+|\nabla_{\sbmx^\prime} f|^2}. 
\end{align*}
%
%
When $\nabla_{\sbmx^\prime} f=\bmzero'$ at some point, $H$ and $|A|^2$ are given by 
\begin{align*}
H=\sum_{i,j\in\Lambda}\delta_{ij}\partial_i\partial_jf=\sum_{i\in\Lambda}\partial_i^2f, \quad
|A|^2=\sum_{i,j,k,\ell\in\Lambda}\delta_{ik}\delta_{j\ell}\partial_i\partial_jf\,\partial_k\partial_\ell f=\sum_{i,j\in\Lambda}(\partial_i\partial_jf)^2. 
\end{align*}
We derive the Laplace-Beltrami operator $\Delta_g$. The definition of the Laplace-Beltrami operator $\Delta_g$ is
\begin{align*}
\Delta_g
=\sum_{i,j\in\Lambda}\frac1{\sqrt{g}}\partial_i\bigl(\sqrt{g}g^{ij}\partial_j\bigr)
=\sum_{i,j\in\Lambda}\left\{g^{ij}\partial_i\partial_j+\frac1{\sqrt{g}}\partial_i\bigl(\sqrt{g}g^{ij}\bigr)\partial_j\right\}.
\end{align*}
Since 
\[
\partial_ig=\partial_i(1+|\nabla_{\sbmx^\prime} f|^2)=2\sum_{k\in\Lambda}\partial_kf\,\partial_i\partial_kf,
\]
it follows that 
\begin{align*}
\sum_{i\in\Lambda}\partial_i(\sqrt{g}g^{ij})
=&\,\sum_{i\in\Lambda}(\partial_i\sqrt{g})g^{ij}+\sum_{i\in\Lambda}\sqrt{g}\,\partial_ig^{ij} \\
=&\,\frac1{2\sqrt{g}}\sum_{i\in\Lambda}(\partial_ig)\left(\delta_{ij}-\frac{\partial_if\,\partial_jf}{g}\right) \\
&\,-\sqrt{g}\sum_{i\in\Lambda}\left(\frac{\partial_i^2f\,\partial_jf+\partial_if\,\partial_i\partial_jf}{g}
-\frac{\partial_if\,\partial_jf\,\partial_ig}{g^2}\right) \\
=&\,\frac1{\sqrt{g}}\sum_{k\in\Lambda}\partial_kf\,\partial_j\partial_kf
+\frac1{g\sqrt{g}}\,\partial_jf\sum_{i,k\in\Lambda}\partial_if\,\partial_kf\,\partial_i\partial_kf \\
&\,-\frac1{\sqrt{g}}\sum_{i\in\Lambda}\bigl(\partial_i^2f\,\partial_jf+\partial_if\,\partial_i\partial_jf\bigr) \\
=&\,\frac1{g\sqrt{g}}\,\partial_jf\sum_{i,k\in\Lambda}\partial_if\,\partial_kf\,\partial_i\partial_kf
-\frac1{\sqrt{g}}\partial_jf\sum_{i\in\Lambda}\partial_i^2f \\
=&\,-\partial_jf\sum_{i,k\in\Lambda}\left(\delta_{ik}-\frac{\partial_if\,\partial_kf}{g}\right)
\frac{\partial_i\partial_kf}{\sqrt{g}} \\
=&\,-H\partial_jf.
\end{align*}
Note that $\nabla_{\sbmx^\prime} f=\bmzero'$ at some point gives $\partial_i(\sqrt{g}g^{ij})=0$. Here we have
\begin{align*}
\partial_k\partial_{\ell}H
=\sum_{i,j\in\Lambda}\bigl\{(\partial_k\partial_{\ell}g^{ij})h_{ij}+\partial_kg^{ij}\partial_{\ell}h_{ij}
+\partial_{\ell}g^{ij}\partial_kh_{ij}+g^{ij}(\partial_k\partial_{\ell}h_{ij})\bigr\}.
\end{align*}
The first partial derivatives of $g^{ij}$ and $h_{ij}$ are given by
\begin{align*}
\partial_kg^{ij}
=&\,-\frac{\partial_i\partial_kf\,\partial_jf+\partial_if\,\partial_j\partial_kf}{g}+\frac{\partial_if\,\partial_jf\,\partial_kg}{g^2}, \\
\partial_kh_{ij}
=&\,\frac{\partial_i\partial_j\partial_kf}{\sqrt{g}}-\frac{\partial_i\partial_jf\,\partial_kg}{2g\sqrt{g}},
\end{align*}
and the second partial derivatives of $g^{ij}$ and $h_{ij}$ are given by
\begin{align*}
\partial_k\partial_{\ell}g^{ij}
=&\,-\frac{\partial_i\partial_k\partial_{\ell}f\,\partial_jf+\partial_i\partial_kf\,\partial_j\partial_{\ell}f
+\partial_i\partial_{\ell}f\,\partial_j\partial_kf+\partial_if\,\partial_j\partial_k\partial_{\ell}f}{g} \\
&\,+\frac{(\partial_i\partial_kf\,\partial_jf+\partial_if\,\partial_j\partial_kf)\partial_{\ell}g}{g^2} \\
&\,+\frac{\partial_i\partial_{\ell}f\,\partial_jf\,\partial_kg+\partial_if\,\partial_j\partial_{\ell}f\,\partial_kg
+\partial_if\,\partial_jf\,\partial_k\partial_{\ell}g}{g^2} \\
&\,-\frac{2\,\partial_if\,\partial_jf\,\partial_kg\,\partial_{\ell}g}{g^3}, \\
\partial_k\partial_{\ell}h_{ij}
=&\,\frac{\partial_i\partial_j\partial_k\partial_{\ell}f}{\sqrt{g}}
-\frac{\partial_i\partial_j\partial_kf\,\partial_{\ell}g}{2g\sqrt{g}}
-\frac{\partial_i\partial_j\partial_{\ell}f\,\partial_kg+\partial_i\partial_jf\,\partial_k\partial_{\ell}g}{2g\sqrt{g}} \\
&\,+\frac{3\,\partial_i\partial_jf\,\partial_kg\,\partial_{\ell}g}{4g^2\sqrt{g}}.
\end{align*}
Note that 
\[
\partial_k\partial_\ell g=2\sum_{m\in\Lambda}\bigl(\partial_k\partial_mf\,\partial_\ell\partial_mf
+\partial_mf\,\partial_k\partial_\ell\partial_mf\bigr).
\]
When $\nabla_{\sbmx^\prime} f=\bmzero'$ at some point, we have
\begin{align*}
&g=1, \quad \partial_kg=0, \quad 
\partial_k\partial_{\ell}g=2\sum_{m\in\Lambda}\partial_k\partial_mf\,\partial_{\ell}\partial_mf, \\
&g_{ij}=g^{ij}=\delta_{ij}, \quad \partial_kg^{ij}=0, \quad
\partial_k\partial_{\ell}g^{ij}=-(\partial_i\partial_kf\,\partial_j\partial_{\ell}f+\partial_i\partial_{\ell}f\,\partial_j\partial_kf),
\end{align*}
so that we obtain 
\begin{align*}
\partial_k\partial_{\ell}H
=&\,-\sum_{i,j\in\Lambda}(\partial_i\partial_kf\,\partial_j\partial_{\ell}f
+\partial_i\partial_{\ell}f\,\partial_j\partial_kf)\partial_i\partial_jf \\
&\,+\sum_{i,j\in\Lambda}\delta_{ij}\left(\partial_i\partial_j\partial_k\partial_{\ell}f
-\partial_i\partial_jf\sum_{m\in\Lambda}\partial_k\partial_mf\,\partial_{\ell}\partial_mf\right).
\end{align*}
This implies 
\begin{align*}
\Delta_gH=&\,\sum_{k,\ell\in\Lambda}g^{k\ell}\partial_k\partial_{\ell}H \\
=&\,\sum_{i,j\in\Lambda}\sum_{k,\ell\in\Lambda}\delta_{ij}\delta_{k\ell}\left(\partial_i\partial_j\partial_k\partial_{\ell}f
-\partial_i\partial_jf\sum_{m\in\Lambda}\partial_k\partial_mf\,\partial_{\ell}\partial_mf\right) \\
&\,-\sum_{i,j\in\Lambda}\sum_{k,\ell\in\Lambda}\delta_{k\ell}(\partial_i\partial_kf\,\partial_j\partial_{\ell}f
+\partial_i\partial_{\ell}f\,\partial_j\partial_kf)\partial_i\partial_jf \\
=&\,\sum_{i,j\in\Lambda}\partial_i^2\partial_j^2f-H|A|^2-2\sum_{i,j,k\in\Lambda}\partial_i\partial_jf\,\partial_i\partial_kf\,\partial_j\partial_kf.
\end{align*}
Consequently, it follows that under the assumption $\nabla_{\sbmx^\prime} f=\bmzero'$ at some point
\begin{align*}
&\hspace*{-10pt}
\Delta_gH+H|A|^2-\frac12H^3 \\[0.05cm]
=&\,\sum_{i,j\in\Lambda}\partial_i^2\partial_j^2f-H|A|^2-2\sum_{i,j,k\in\Lambda}\partial_i\partial_jf\,\partial_i\partial_kf\,\partial_j\partial_kf
+H|A|^2-\frac12H^3 \\
=&\,\sum_{i,j\in\Lambda}\partial_i^2\partial_j^2f-\frac12H^3-2\sum_{i,j,k\in\Lambda}\partial_i\partial_jf\,\partial_i\partial_kf\,\partial_j\partial_kf.
\end{align*}
Since $A=(\partial_i\partial_jf)$ is a symmetric matrix, there exists an orthogonal matrix $P_{N-1}=(p_{ij})$ such that
\[
P_{N-1}^{-1}AP_{N-1}=\left(\begin{array}{ccc} \kappa_1& & \\ &\ddots& \\ & &\kappa_{N-1} \end{array}\right),
\]
where $\kappa_i\,(i=1,\cdots,N-1)$ are the principle curvatures. Taking account of $P_{N-1}^{-1}={\,}^tP_{N-1}$, it 
follows that
\[
A=P_{N-1}\left(\begin{array}{ccc} \kappa_1& & \\ &\ddots& \\ & &\kappa_{N-1} \end{array}\right)P_{N-1}^{-1}
=P_{N-1}\left(\begin{array}{ccc} \kappa_1& & \\ &\ddots& \\ & &\kappa_{N-1} \end{array}\right){\,}^tP_{N-1},
\]
so that we have
\[
\partial_i\partial_jf=\sum_{m\in\Lambda}\kappa_mp_{im}p_{jm}.
\]
Since $\sum\limits_kp_{ki}p_{kj}=\delta_{ij}$, this implies 
\begin{align*}
&\sum_{i,j,k\in\Lambda}\partial_i\partial_jf\,\partial_i\partial_kf\,\partial_j\partial_kf \\
&=\sum_{i,j,k\in\Lambda}\sum_{m_1,m_2,m_3\in\Lambda}\kappa_{m_1}\kappa_{m_2}\kappa_{m_3}p_{i,m_1}p_{j,m_1}p_{i,m_2}p_{k,m_2}p_{j,m_3}p_{k,m_3} \\
&=\sum_{m_1,m_2,m_3\in\Lambda}\kappa_{m_1}\kappa_{m_2}\kappa_{m_3}\left(\sum_{i\in\Lambda}p_{i,m_1}p_{i,m_2}\right)
\left(\sum_{j\in\Lambda}p_{j,m_1}p_{j,m_3}\right)
\left(\sum_{k\in\Lambda}p_{k,m_2}p_{k,m_3}\right) \\
&=\sum_{m_1,m_2,m_3\in\Lambda}\kappa_{m_1}\kappa_{m_2}\kappa_{m_3}\delta_{m_1m_2}\delta_{m_1m_3}\delta_{m_2m_3} \\
&=\sum_{m\in\Lambda}\kappa_m^3.
\end{align*}
Thus we are led to
\begin{align*}
&\hspace*{-10pt}
\Delta_gH+H|A|^2-\frac12H^3 \\
=&\,\sum_{i,j\in\Lambda}\partial_i^2\partial_j^2f-\frac12\left(\sum_{m\in\Lambda}\kappa_m\right)^3
-2\sum_{m\in\Lambda}\kappa_m^3 \\
=&\,\sum_{i\in\Lambda}\partial_i^4f+2\sum_{\substack{i,j\in\Lambda \\ i<j}}\partial_i^2\partial_j^2f \\
&\,-\frac52\sum_{m\in\Lambda}\kappa_m^3
-\frac32\sum_{\substack{m_1,m_2\in\Lambda \\ m_1\ne m_2}}\kappa_{m_1}^2\kappa_{m_2}
-3\sum_{\substack{m_1,m_2,m_3\in\Lambda \\ m_1<m_2<m_3}}\kappa_{m_1}\kappa_{m_2}\kappa_{m_3}.
\end{align*}

\section{Proof of Theorem \ref{thm:Cui}.}

Theorem \ref{thm:Cui} with $ \lambda=0 $ is the same as Cui\,\cite[Theorem 3.2]{cui;01}.  However, in \cite{cui;01} the author does 
not state the dependence of the constant $ C $ on $ \al\in\bZ^N $ and $ m\in\bN $. 
Noting this point, we give a proof of Theorem \ref{thm:Cui} in this subsection.  

First, we prepare some notations, according to \cite{cui;01}.  Set 
\[
P_{2k}(\textbf{i}\bmzeta):=\{(\textbf{i}\zeta_1)^2+\cdots+(\textbf{i}\zeta_N)^2\}^k
=(-1)^k(\zeta_1^2+\cdots+\zeta_N^2)^k
\] 
for $\bmzeta=(\zeta_1,\cdots,\zeta_N)\in\bC^N$ and $k\in\mathbb{Z}_+$.  Since it follows from direct calculations 
that there exist $\delta\in(0,1)$, $K_1>1$ and $K_2>1$ such that 
\[
\text{Re}\,\{-P_4(\textbf{i}\bmzeta)-\lambda P_2(\textbf{i}\bmzeta)\}
\leq -\delta|\text{Re}\,\bmzeta|^4+K_1|\text{Im}\,\bmzeta|^4+K_2|\lambda|^2
\]
for $\lambda\in\mathbb{R}$, where $\text{Re}\,\bmzeta:=(\text{Re}\,\zeta_1,\cdots,\text{Re}\,\zeta_N) $ and 
$ \text{Im}\,\bmzeta:=(\text{Im}\,\zeta_1,\cdots,\text{Im}\,\zeta_N) $, we are led to
\begin{equation}
\label{ineqap0201}
|e^{-P_4(\textbf{i}\sbmzeta)-\lambda P_2(\textbf{i}\sbmzeta)}|
\leq e^{-\delta|\text{Re}\,\sbmzeta|^4+K_1|\text{Im}\,\sbmzeta|^4+K_2|\lambda|^2} \quad
\mbox{for all}\,\ \bmzeta\in\bC^N.    
\end{equation}
Note that 
\[
G_N(\bmx,t)=c_N\int_{\bR^N}e^{-(P_4(\textbf{i}\sbmxi)+\lambda P_2(\textbf{i}\sbmxi))t+\textbf{i}\langle \sbmx,\sbmxi\rangle_N}\,d\bmxi.
\] 
We are now in a position to the proof of Theorem \ref{thm:Cui}. 

\medskip
\noindent
{\bf Proof of Theorem \ref{thm:Cui}.} 
The proof follows from the argument in that of \cite[Theorem 3.2]{cui;01}.    
It follows from Cauchy's integral theorem that  
\begin{align*}
&D_{\sbmx}^\al(-\Delta_{\sbmx})^mG_N(\bmx,t) \\
&=c_N\int_{\bR^N}(\textbf{i}\bmxi-\bmeta)^{\al}\biggl\{-\sum_{j=1}^N(\textbf{i}\xi_j-\eta_j)^2\biggr\}^m
e^{-(P_4(\sbmxi+\textbf{i}\sbmeta)+\lambda P_2(\sbmxi+\textbf{i}\sbmeta))t
+\textbf{i}\langle\sbmx,\sbmxi\rangle_N-\langle\sbmx,\sbmeta\rangle_N}\,d\bmxi, 
\end{align*}
where $\bmeta\in\bR^N$ is arbitrary and independent of $\bmxi\in\bR^N$.  Then \eqref{ineqap0201} implies that 
\begin{align*}
&|D_{\sbmx}^\al(-\Delta_{\sbmx})^mG_N(\bmx,t)| \\
&\leq c_Ne^{-\langle\sbmx,\sbmeta\rangle_N}\int_{\bR^N}(|\bmxi|+|\bmeta|)^{|\al|+2m}
|e^{-(P_4(\sbmxi+\textbf{i}\sbmeta)+\lambda P_2(\sbmxi+\textbf{i}\sbmeta))t}|\,d\bmxi \\
&\leq c_Ne^{-\langle\sbmx,\sbmeta\rangle_N+K_1|\sbmeta|^4t+K_2|\lambda|^2t}
\int_{\bR^N}(|\bmxi|+|\bmeta|)^{|\al|+2m}e^{-\delta|\sbmxi|^4t}d\bmxi \\
&\leq c_Ne^{-\langle\sbmx,\sbmeta\rangle_N+K_1|\sbmeta|^4t+K_2|\lambda|^2t}
\sum_{k=0}^{|\al|+2m}\binom{|\alpha|+2m}{k}|\eta|^{|\al|+2m-k}\int_{\bR^N}|\bmxi|^{k}e^{-\delta|\sbmxi|^4t}\,d\bmxi.  
\end{align*}
Applying the change of variable on the polar coordinate, we have
\begin{align*}
\int_{\bR^N}|\bmxi|^{k}e^{-\delta|\sbmxi|^4t}\,d\bmxi
=&\,\omega_{N-1}\int_0^\infty r^{k+N-1}e^{-\delta r^4t}\,dr \\
=&\,\frac{\omega_{N-1}}4\Gamma\Bigl(\frac{k+N}4\Bigr)(\delta t)^{-(N+k)/4} \\
\le&\,\frac{\omega_{N-1}}4\max\biggl\{\Gamma\Bigl(\frac{|\alpha|+2m+N}4\Bigr),\Gamma\Bigl(\frac{N}4\Bigr)\biggr\}(\delta t)^{-(N+k)/4},
\end{align*}
where $\omega_{N-1}$ is the area of the $(N-1)$-dimensional unit ball. Taking account of the fact that the minimum of 
$\Gamma(s)$ exists and is positive for $s>0$, we see that
\begin{align*}
&|D_{\sbmx}^\al(-\Delta_{\sbmx})^mG_N(\bmx,t)| \\
&\le C\,\Gamma\Bigl(\frac{|\alpha|+2m+N}4\Bigr)t^{-N/4}e^{-\langle\sbmx,\sbmeta\rangle_N+K_1|\sbmeta|^4t+K_2|\lambda|^2t}
\bigl\{(\delta t)^{-1/4}+|\bmeta|\bigr\}^{|\alpha|+2m},
\end{align*}
where $C>0$ is a constant depending only on $N$. Since $\bmeta\in\bR^N$ is arbitrary, we can choose 
$\bmeta:=(4K_1)^{-1/3}t^{-1/3}|\bmx|^{-2/3}\bmx$.  Then it follows that 
\begin{align*}
&e^{-\langle\sbmx,\sbmeta\rangle_N+K_1|\sbmeta|^4t+K_2|\lambda|^2t}\bigl\{(\delta t)^{-1/4}+|\bmeta|\bigr\}^{|\alpha|+2m} \\
&=e^{-\mu(|x|^4/t)^{1/3}+K_2|\lambda|^2t}
\biggl\{\frac1{(\delta t)^{1/4}}+\frac1{(4K_1)^{1/3}}\biggl(\frac{|\bmx|}{t}\biggr)^{1/3}\biggr\}^{|\alpha|+2m} \\
&=\nu^{(|\alpha|+2m)/4}\,t^{-(|\alpha|+2m)/4}e^{-\mu(|x|^4/t)^{1/3}+K_2|\lambda|^2t}
\biggl\{1+\frac{\delta^{1/4}}{(4K_1)^{1/3}}\biggl(\frac{|\bmx|}{t^{1/4}}\biggr)^{1/3}\biggr\}^{|\alpha|+2m} \\
&\le\nu^{(|\alpha|+2m)/4}\,t^{-(|\alpha|+2m)/4}e^{-\mu(|x|^4/t)^{1/3}+K_2|\lambda|^2t}
\biggl\{1+\biggl(\frac{|\bmx|}{t^{1/4}}\biggr)^{1/3}\biggr\}^{|\alpha|+2m},
\end{align*}
where 
\[
\mu=\frac34(4K_1)^{-1/3}\in(0,1), \quad \nu=\frac1{\delta}>1.
\]
In the last inequality, we have used $\delta\in(0,1)$ and $K_1>1$. Since $1+r^{1/3}\le2^{2/3}(1+r)^{1/3}$ 
for $r>0$, we are led to the desired result. \qed



\bigskip\noindent
(Katsuyuki Ishii) {\textsc{Graduate School of Maritime Sciences, Kobe University, 
5-1-1 Fukaeminami-machi, Higashinada-ku, Kobe 658-0022, Japan}} \\
{\it Email address}: {\texttt{ishii@maritime.kobe-u.ac.jp}} \\[0.325cm]
(Yoshihito Kohsaka) {\textsc{Graduate School of Maritime Sciences, Kobe University, 
5-1-1 Fukaeminami-machi, Higashinada-ku, Kobe 658-0022, Japan}} \\
{\it Email address}: {\texttt{kohsaka@maritime.kobe-u.ac.jp}} \\[0.325cm]
(Nobuhito Miyake) {\textsc{Faculty of Mathematics, Kyushu University, 
744, Motooka, Nishi-ku, Fukuoka-shi, Fukuoka 819-0395, Japan}} \\
{\it Email address}: {\texttt{miyake@math.kyushu-u.ac.jp}} \\[0.325cm]
(Koya Sakakibara) {\textsc{Faculty of Mathematics and Physics, Institute of Science and Engineering, Kanazawa University, 
Kakuma-machi, Kanazawa 920-1192, Japan; 
RIKEN Interdisciplinary Theoretical and Mathematical Sciences Program (iTHEMS), 
2-1 Hirosawa, Wako, Saitama 351-0198 Japan}} \\
{\it Email address}: {\texttt{ksakaki@se.kanazawa-u.ac.jp}}

\end{document}